\documentclass{amsart}
\newcommand{\N}{{\Bbb N}}

\newcommand{\Z}{{\Bbb Z}}

\newcommand{\C}{{\Bbb C}}

\renewcommand{\P}{{\Bbb P}}

\newcommand{\kf}{{\mathcal F}}
\newcommand{\kg}{{\mathcal G}}

\newcommand{\kj}{{\mathcal J}}

\newcommand{\kn}{{\mathcal N}}
\newcommand{\ko}{{\mathcal O}}

\newcommand{\ks}{{\mathcal S}}
\newcommand{\kt}{{\mathcal T}}

\newcommand{\Virr}{\mbox{$V_{irr}(d;\,d_1,\ldots,d_r;\,k)$}}
\newcommand{\tilN}{\tilde{\kn}_{\tilde{H}/\P^2_r}^{\Sigma_0}}
\newtheorem{lemma}[equation]{Lemma}
\newtheorem{proposition}[equation]{Proposition}
\newtheorem{theorem}{Theorem}
\newtheorem{corollary}[equation]{Corollary}

\theoremstyle{definition}

\theoremstyle{remark}
\newtheorem{remark}[equation]{Remark}

\input amssym.def
\input amssym
\begin{document}

\title[Families of Nodal Curves]{Geometry of Families of Nodal Curves on
 the blown--up projective plane}
\author[G.-M. Greuel]{Gert-Martin Greuel}
\address{
Universit\"at Kaiserslautern\\
Fachbereich Mathematik\\
Erwin-Schr\"odinger-Stra\ss e\\
D -- 67663 Kaiserslautern}
\email{greuel@mathematik.uni-kl.de}
\author[C. Lossen]{Christoph Lossen}
\address{
Universit\"at Kaiserslautern\\
Fachbereich Mathematik\\
Erwin-Schr\"odinger-Stra\ss e\\
D -- 67663 Kaiserslautern}
\email{lossen@mathematik.uni-kl.de}
\author[E. Shustin]{Eugenii Shustin}
\address{
Tel Aviv University\\
School of Mathematical Sciences\\
Ramat Aviv\\
ISR -- Tel Aviv 69978}
\email{shustin@math.tau.ac.il}

\begin{abstract}
Let  $\P^2_r\,$ be the projective plane blown up at $r$
generic points. De\-note by \mbox{$E_0,E_1,\ldots ,E_r$} the strict transform
of a generic straight line on $\P^2$ and the exceptional divisors of the
blown--up points on $\P^2_r$ respectively. We consider the variety \Virr\ of
all irreducible curves $C$ in \mbox{$|dE_0-\sum_{i=1}^{r}
  d_iE_i|$} with $k$ nodes as the only singularities and give
asymp\-to\-ti\-cally nearly optimal sufficient conditions for its smoothness,
irreducibility and non--emptyness. Moreover,
we extend our conditions for the smooth\-ness and the irreducibility on
families of
reducible curves. For \mbox{$r\leq 9$} we give the
complete answer concerning the existence of nodal curves in \Virr .
\end{abstract}

\maketitle

\section*{Introduction}

We deal with the following general problem: given a
smooth rational surface $S$ and a divisor $D$ on $S$, when is the variety
$V_{irr}(D,k)$ of nodal irreducible curves in the complete linear system $|D|$
with a fixed number $k$ of nodes non--empty, when non--singular and when
irreducible? For \mbox{$S=\P^2$}, these questions are completely answered by
the classical result of F.~Severi (\cite{Sev}), stating that the variety
$V_{irr}(dH,k)$ of
irreducible curves of degree $d$ having $k$ nodes is non--empty and smooth
exactly if $$0\leq k \leq \frac{(d-1)(d-2)}{2}\:,$$
and the result of J.~Harris (\cite{Har}), stating that $V_{irr}(dH,k)$ is
always irreducible.\\
A modification of Severi's method did lead to a sufficient
(smoothness--)criterion for general smooth rational surfaces $S$
(\cite{Ta1,Nob}): let \mbox{$C_0\subset S$} be a smooth
irreducible curve, let \mbox{$C\in |C_0|$} be a reduced (nodal) curve with
precisely $k$ nodes, such that \mbox{$C=C_1\cup \ldots \cup C_s$}, $C_i$
irreducible and
\begin{equation}
  \label{0.1}
  K_S\cdot C_i < 0
\end{equation}
for each \mbox{$1\leq i\leq s$}, then the variety $V_{irr}(C_0,k)$ of
irreducible curves in $|C_0|$ having precisely $k$ nodes is smooth (see
\cite{Ta2,GrM,GrK,GrL} for generalizations to other surfaces). Moreover,
in those cases each node of $C$ can be smoothed independently.\\
In this paper, we concentrate on the case \mbox{$S=\P^2_r$}, the projective
plane blown up at $r$ {\em generic} points \mbox{$p_1,\ldots,p_r$}. Let
\mbox{$E_0,E_1,\ldots ,E_r$} denote the strict transform of a generic
straight line on $\P^2$ and the exceptional divisors of the blown--up points on
$\P^2_r$, respectively. Then for an irreducible nodal curve
\mbox{$C\in|dE_0-\sum_{i=1}^{r} d_iE_i|\,$} the condition (\ref{0.1}) reads as
\begin{equation}
  \label{0.2}
  3d > \sum_{i=1}^r d_i\:.
\end{equation}
In the blown--down situation, such a curve $C$ corresponds to a plane curve of
degree $d$ having (not necessarily ordinary) $d_i$--fold points at $p_i$,
\mbox{$1\leq i\leq r$}, and \mbox{$k'\leq k$} nodes outside. For the variety of
irreducible plane curves of fixed topological (or analytic) type, E.~Shustin
gives in (\cite{Sh2}) an asymptotically improved sufficient condition for
smoothness and irreducibility:
$$
  \alpha d^2+o(d^2) > \sum \sigma(S_i)\:,
$$
where $\sigma$ denotes some positive invariant of the singular points. In our
case, \mbox{keeping} $k$, $d$ and the $d_i$ \mbox{$(1\leq i \leq r)$} fixed,
certainly $k'$ and
the topological types of the multiple points may vary. Nevertheless, we shall
obtain sufficient conditions for the smoothness and the
irreducibility of the same type, that is, with the same
exponent in $d$. Moreover, we should like to emphasize that we can extend them
on families of reducible curves
- --- in contrast to A.~Tannenbaum's result for K3-surfaces (cf.~\cite{Ta2}).\\
In section 3, we shall give a complete answer for the existence problem in case
of \mbox{$r\leq 9$} blown--up points (Theorems \ref{4.1A} and \ref{4.1B}). For
\mbox{$r\geq 10$}, we obtain an {\em exponentially optimal} sufficient
condition (Corollary \ref{3.1.4}), that is, of the same exponent in
$d$ as the known restrictions for the existence of the corresponding plane
curves with $d_i$--fold singularities $S_i$ \mbox{$(1\leq i \leq r)$} and $k'$
nodes \mbox{$S_{r+1},\ldots,S_{r+k'}$} (from
Pl\"ucker formulae to inequalities by Varchenko \cite{Var} and Ivinskis
\cite{Ivi,HiF}). These restrictions are of type
$$
  \alpha_2 d^2+\alpha_1 d+\alpha_0 > \sum_{i=1}^{r+k'}
  \sigma(S_i)\;\;\;\;\;\;\;\;\;\; (\alpha_2=\mbox{const}>0)
$$
with $\sigma$ some positive invariant depending, at most, quadratically on $d$.
Our result improves for the given situation the only previously known
(general) existence criterion (in \cite{Sh1}):
$$ \frac{(d+3)^2}{2} \geq \sum_{i=1}^{r+k'} (\mu(S_i)+4)(\mu(S_i)+5)\:, $$
which is not exponentially optimal since the right--hand side may be of order
four in $d$.
For the proof, we combine a modification of the method of
A.~Hirschowitz in \cite{Hir} and the smoothing of nodes (cf.~\cite{Ta1}).

\section*{Notation and terminology}

Throughout this article we consider all objects to be defined over an
algebraically closed field {\bf $K$} of characteristic zero. We use the
following notations:
\begin{itemize}
\itemsep0.1cm
\item $\P^2_r$ --- the projective plane blown up at r generic
points \mbox{$p_1,\ldots,p_r$}.
\item $\frak{m}_{z_\nu}$ --- the maximal ideal in the local ring
  $\ko_{\P^2_r,z_{\nu} }$, \mbox{$z_{\nu}\in \P^2_r$}.
\item $E_0$ --- the strict transform of a generic straight line (in $\P^2$).
\item $E_i$ \mbox{$(1\leq i\leq r)$} --- the exceptional divisor of the
  blown--up  point $p_i$  on $\P^2_r$.
\item \Virr\ ---
the variety of all irreducible curves $C$ in the linear system
\mbox{$|dE_0-\sum_{i=1}^{r} d_iE_i|\,$}
having k nodes as their only singularities.
\end{itemize}
Furthermore, for a reduced nodal curve \mbox{$C\subset \P^2_r$},
\mbox{$\,C=C_1\cup\ldots\cup C_s$} ($C_i$ irreducible), having precisely $k$
nodes as their only singularities and a divisor $D$ on $\P^2_r$, we denote:
\begin{itemize}
\item \mbox{$V(|D|;C)$} --- the variety of all reduced curves
  \mbox{$\tilde{C} =\tilde{C}_1 \cup\ldots\cup \tilde{C}_s $} in the linear
  system
  $|D|$ with precisely $k$ nodes as their only singularities, whose components
  $\tilde{C}_i $ have the same type (that is, are in the same linear system and
  have the same number of nodes) as the components $C_i$ (\mbox{$1\leq i\leq
  s$}) of $C$.
\end{itemize}

\section{Smoothness}

\subsection{Formulation of the result}
\setcounter{equation}{0}
For $\P^2_r\,$, \mbox{$r\leq 8\,$}, condition (\ref{0.1}) is fulfilled for each
irreducible curve $C$, hence \Virr\ is always smooth. In case $r=9$ and
\mbox{$C\in V_{irr}(d;\,d_1,\ldots,d_9;\,k)$}, (\ref{0.1}) reads
 $$ 3d > \sum_{i=1}^9 d_i\:, $$
which is satisfied exactly if $C$ is not the (unique) smooth cubic through
\mbox{$p_1,\ldots, p_9$}. Thus \mbox{$V_{irr}(d;\,d_1,\ldots,d_9;\,k)$} is
always smooth, too.
In this section we shall prove:

\begin{theorem}[Smoothness Theorem] \label{1.1.1}
Let \mbox{$r\geq 10$}, and let the positive integers
  \mbox{$d;\,d_1,\ldots,d_r$} satisfy the two (smoothness) conditions
\renewcommand{\arraystretch}{1.5}
\begin{eqnarray}
    \left[ \; \sqrt{2k}\;\, \right] & < & \frac{d}{2} +3-\frac{\sqrt{2} }{2}
    \sqrt{\sum_{i=1}^r (d_i+2)^2} \label{Bed1} \\
    \left[ \; \sqrt{2k}\;\, \right] & < & d+3-\sqrt{2}\sqrt{2+\sum_{i=1}^r
      (d_i+2)(d_i+1)} \:\:.
    \label{Bed2}
  \end{eqnarray}
Then \Virr\  is smooth and has the T--property (that is, each germ of \Virr\ is
a transversal intersection of germs of equisingular strata corresponding to the
$k$ nodes).

Let \mbox{$C \subset \P^2_r$} be a reduced nodal curve, \mbox{$C\in
    |D|$}. If, for each irreducible component \mbox{$C_{\nu} \in
    |d^{(\nu)}E_0-\sum_{i=1}^{r} d^{(\nu)}_iE_i|\,$} of $C$ (having precisely
$k^{(\nu)}$ nodes), the two smoothness conditions are fulfilled, then
  \mbox{$V(|D|;C)$} is smooth and has the T--property.
\end{theorem}

\subsection{Vanishing criteria} \label{1.2}
\setcounter{equation}{0}

We introduce the vanishing criteria which we
shall mainly use
in the proof of the Smoothness and the Irreducibility Theorem (in the next
paragraph):\\
For a curve \mbox{$C\in \,$\Virr} and a subset \mbox{$\Sigma_0 \subset$
  Sing$\,C$} we define
the sheaf $\kt ^1_{C,\Sigma_0}$ to be the skyscraper sheaf concentrated at
\mbox{$\mbox{Sing}\,C=\{z_1,\ldots ,z_k\}$} with stalks
\renewcommand{\arraystretch}{1.3}
$$\left( \kt ^1_{C,\Sigma_0} \right) _{z_\nu} := \Bigg\{
    \begin{array}{ccl}
      \ko_{\P^2_r,z_\nu}/\frak{m}_{z_\nu} & \mbox{for }  & z_\nu \in
      \mbox{Sing}\,C-\Sigma_0 \\
      \ko_{\P^2_r,z_\nu}/\frak{m}_{z_\nu}^2 & \mbox{for } & z_\nu \in
      \Sigma_0 \:\:.\\
    \end{array}
$$
Furthermore, put
$$ \kn_{C/\P^2_r}^{\Sigma_0} :=\mbox{Ker} \left(\kn_{C/\P^2_r}\longrightarrow
  \kt^1_{C,\Sigma_0} \right) .$$

\renewcommand{\labelenumi}{(\Alph{enumi})}
\begin{proposition} \label{1.2.1}
  Let \mbox{$C\in |dE_0-\sum_{i=1}^r d_iE_i|$} be a reduced nodal curve having
  $k$ nodes as its only singularities,
  \mbox{$\tilde{C}\sim dE_0-\sum_{i=1}^{r} \tilde{d_i}E_i\,$},
  where \mbox{$\tilde{d_i}\geq d_i$} for \mbox{$1\leq i\leq r$}, and
  \mbox{$\Sigma_0 \subset$
  Sing$\,C=\{z_1,\ldots z_k\}$}. Moreover, let $\tilde{H}$ be a reduced
  curve whose local equations map to $0\in\left( \kt ^1_{C,\Sigma_0} \right)
  _{z_\nu}$ for $1\leq \nu \leq k$.\\
  Then $H^1(C,\kn_{C/\P^2_r}^{\Sigma_0})$ vanishes, if the following
conditions are satisfied:
  \begin{enumerate}
  \itemsep0.2cm
  \item \label{A} $\;H^1(\P^2_r,\ko_{\P^2_r}(C))=0$
  \item \label{B} $\;H^1(\P^2_r,\ko_{\P^2_r}(\tilde{C}-\tilde{H}))=0$
  \item \label{C}
    $\;H^1(\tilde{H},\tilde{\kn}_{\tilde{H}/\P^2_r}^{\Sigma_0})=0$
  \end{enumerate}
  where $\;\tilde{\kn}_{\tilde{H}/\P^2_r}^{\Sigma_0}:= \mbox{Ker}
  \left(\ko_{\P^2_r}(\tilde{C})\otimes \ko_{\tilde{H}} \longrightarrow
  \kt^1_{C,\Sigma_0} \right)\,$.
\end{proposition}

\begin{proof}
We have an exact sequence
$$  \ldots \rightarrow H^1(\P^2_r,\ko_{\P^2_r}(C))\rightarrow
H^1(C,\kn_{C/\P^2_r} )\rightarrow H^2(\P^2_r,\ko_{\P^2_r})\rightarrow
\ldots
$$
where, by Serre duality,
$H^2(\P^2_r,\ko_{\P^2_r})=H^0(\P^2_r,K_{\P^2_r})=0\,$. Hence, the statement of
the proposition follows immediately from the following commutative diagram with
exact columns and diagonal

\unitlength0.7cm
\begin{picture}(18,8)
\put(0.3,7.0){\makebox(3.2,0.8){$\vdots$}}
\put(3.0,6.8){\makebox(7.0,0.8){$H^0(\ko_{\P^2_r}(C))\supset
  H^0(\ko_{\P^2_r}(\tilde{C}))$}}
\put(9.7,7.0){\makebox(3.2,0.8){$\vdots$}}
\put(1.9,6.7){\vector(0,-1){0.6}}
\put(4.1,6.7){\vector(-2,-1){1.25}}
\put(9.0,6.7){\vector(2,-1){1.25}}
\put(11.3,6.7){\vector(0,-1){0.6}}
\put(0.3,5.2){\makebox(3.2,0.8){$H^0(\kn_{C/\P^2_r})$}}
\put(9.7,5.2){\makebox(3.2,0.8){$H^0(\ko_{\P^2_r}(\tilde{C})\!\otimes\!
\ko_{\tilde{H}})$}}
\put(1.9,5.1){\vector(0,-1){0.6}}
\put(11.3,5.1){\vector(0,-1){0.6}}
\put(12.2,5.1){\vector(2,-1){1.25}}
\put(0.3,3.6){\makebox(3.2,0.8){$H^0(\kt^1_{C,\Sigma_0})$}}
\put(3.6,3.94){\line(1,0){6.1}}
\put(3.6,4.06){\line(1,0){6.1}}
\put(9.7,3.6){\makebox(3.2,0.8){$H^0(\kt^1_{C,\Sigma_0})$}}
\put(13.8,3.6){\makebox(3.2,0.8){$H^1(\ko_{\P^2_r}(\tilde{C}\!
    -\!\tilde{H}))$}}
\put(15.47,3.4){\line(0,-1){0.3}}
\put(15.33,3.4){\line(0,-1){0.3}}
\put(13.8,2.3){\makebox(3.2,0.8){$0$}}
\put(1.9,3.5){\vector(0,-1){0.6}}
\put(11.3,3.5){\vector(0,-1){0.6}}
\put(0.3,1.9){\makebox(3.2,0.8){$H^1(\kn_{C/\P^2_r}^{\Sigma_0})$}}
\put(9.7,1.9){\makebox(3.2,0.8){$
    H^1(\tilde{\kn}_{\tilde{H}/\P^2_r}^{\Sigma_0})$}}
\put(1.9,1.8){\vector(0,-1){0.6}}
\put(11.37,1.8){\line(0,-1){0.3}}
\put(11.23,1.8){\line(0,-1){0.3}}
\put(9.7,0.7){\makebox(3.2,0.8){$0$}}
\put(0.3,0.3){\makebox(3.2,0.8){$H^1(\kn_{C/\P^2_r})$}}
\put(3.6,0.3){\makebox(0.6,0.8){$=\;0$}}
\put(16.0,0.0){\makebox(1.0,0.8){\qed}}
\end{picture}
\renewcommand{\qed}{}\end{proof}

In the following we shall obtain the vanishing properties (A) -- (C) by
applying two well--known criteria:

\begin{proposition}[Hirschowitz-Criterion, \cite{Hir}] \label{Hirschowitz}
  Let \mbox{$C\sim dE_0-\sum_{i=1}^r d_iE_i$}, where $d\,,d_1,\ldots d_r$ are
  non--negative integers satisfying
  \begin{eqnarray} \label{1.2.3}
    \sum_{i=1}^r \frac{d_i(d_i+1)}{2} & < & \left[\frac{(d+3)^2}{4}\right]\:,
  \end{eqnarray}
  then $H^1(\P^2_r,\ko_{\P^2_r} (C))=0$.

  In a more special situation, let $S_r^I$ be the projective plane blown up at
  r points \mbox{$p_1,\ldots ,p_r$} where $p_i$, \mbox{$i\in I$}, lie on a
  line, and all the
  other points are in generic position. Let $C$ be as above such that
  \begin{eqnarray} \label{1.2.3b}
    \sum_{i\in I} d_i & \leq & d+1
  \end{eqnarray}
  and condition $($\ref{1.2.3}$)$ holds, then \mbox{$H^1(S_r^I,\ko_{S^I_r}
    (C))=0$}.
\end{proposition}

\begin{proposition}[\cite{GrK}] \label{Greuel}  Let $S$
  be a smooth surface,
  \mbox{$C\subset S$} a compact reduced curve, $\kf $ a torsionfree coherent
  $\ko_C$-module which has rank 1 on each irreducible component $C_i$ of $C$
  \mbox{$(1\leq i\leq s)$}. Then \mbox{$H^1(C,\kf)=0$} if for \mbox{$1\leq
    i\leq s$}
 \begin{equation}
    \label{1.2.5}
    \chi (\overline{\kf \otimes \ko_{C_i}})  >  \chi (\omega_C \otimes
    \ko_{C_i}) - \mbox{isod}_{C_i}(\kf, \ko_C)\:.
  \end{equation}
  Here $\overline{\raisebox{0.2cm}{\ \ \ } }$ denotes reduction modulo torsion,
  $\omega_C $ the
  dualizing sheaf and the isomorphism defect $\mbox{isod}_{C_i}(\kf, \ko_C)$ is
  defined to be the sum of all
$$ \:\mbox{isod}_{C_i,x}(\kf, \ko_C):=\mbox{min}\,(\mbox{dim}_{\,\C}\:
\mbox{coker}\,(\varphi
_{C_i}:\,(\overline{\kf \otimes \ko_{C_i}})_x \rightarrow \ko_{C_i,x}))\,,$$
$x\in C_i\,$, where the minimum is taken over all $\varphi_{C_i}$ which are
induced by local homomorphisms \mbox{$\,\varphi :\,\kf_x\rightarrow
  \ko_{C,x}$}.
\end{proposition}

\begin{corollary} \label{1.2.6}
  If, in the situation of Proposition \ref{1.2.1}, we have for each irreducible
  component $\tilde{H_i}\,$ \mbox{$(1\leq i\leq s)$} of $\tilde{H} $
  \begin{eqnarray}
    \label{1.2.7}
    && \tilde{H}_i\cdot (\tilde{C}\!-\!\tilde{H}\!-\!K_{\P^2_r})\; >\;
    \#\,(\mbox{Sing}\,C\cap \tilde{H}_i)\;:= \! \sum_{z\in
    \mbox{\footnotesize Sing}\,C} \!\mbox{multiplicity}\,(\tilde{H}_i,z)
  \end{eqnarray}
  then $H^1(\tilde{H},\tilde{\kn}_{\tilde{H}/\P^2_r}^{\Sigma_0})\,$ vanishes.
\end{corollary}

\begin{proof}
Applying the Riemann--Roch--Theorem and the adjunction formula,
condition (\ref{1.2.5}) reads
$$\mbox{deg} \,(\overline{\tilN \otimes
  \ko_{\tilde{H}_i}})\;>\;(K_{\P^2_r}+\tilde{H})\cdot \tilde{H}_i -
  \mbox{isod}_{\tilde{H}_i}(\tilN, \ko_{\tilde{H}})\:. $$
The exact sequence \mbox{$\;0\rightarrow \tilN \rightarrow
  \ko_{\P^2_r}(\tilde{C}) \otimes
\ko_{\tilde{H}} \rightarrow \kt^1_{C,\Sigma_0} \rightarrow 0\;$} implies
$$ \mbox{deg}\,(\overline{\tilN \otimes \ko_{\tilde{H}_i}})\:=\:
\mbox{deg}\,(\ko_{\P^2_r}(\tilde{C}) \otimes \ko_{\tilde{H}_i}) -
\chi\,(\kt^1_{C,\Sigma_0} \otimes \ko_{\tilde{H}_i})\:. $$
Finally, an easy consideration shows that
\begin{eqnarray*}
\lefteqn{\chi\,(\kt^1_{C,\Sigma_0} \otimes
  \ko_{\tilde{H}_i})-\mbox{isod}_{\tilde{H}_i}(\tilN,
\ko_{\tilde{H}}) } \hspace{1.0cm}\\
& & = \sum_{z\in \,\mbox{\footnotesize Sing}\,C}
dim_{\C}(\kt^1_{C,\Sigma_0}
\otimes \ko_{\tilde{H}_i})_z-isod_{C_i}(\tilN,\ko_{\tilde{H}})\\
& & \leq \#\,(\mbox{Sing}\:C\:\cap \: \tilde{H}_i)\:. \qed
\end{eqnarray*}
\renewcommand{\qed}{}\end{proof}

\subsection{Proof of the Smoothness Theorem}
\setcounter{equation}{0}

Following (\cite{GrK}, Theorem 6.1), it is sufficient
to show that the first
cohomology group $H^1(C,\kn_{C/\P^2_r}^{\;\emptyset})$ of the sheaf
$$ \kn_{C/\P^2_r}^{\;\emptyset}=\mbox{Ker}\,(\kn_{C/\P^2_r}\longrightarrow
\kt^1_{C} )$$
vanishes, where $\kt^1_{C}$ denotes the skyscraper sheaf concentrated
in the singular set \mbox{Sing$\,C=\{z_1,\ldots,z_k\}$} with stalk in
$z_\nu$
$$(\kt^1_{C})_{z_{\nu}} =
\ko_{\P^2_r,z_{\nu}}/\frak{m}_{z_{\nu}}\:.$$
If the reduced curve \mbox{$C\subset \P^2_r$} decomposes as \mbox{$C=C'\cup
  C''$}, then we can consider the exact sequence
$$ 0\longrightarrow \kn_{C'/\P^2_r}\oplus \kn_{C''/\P^2_r}
\stackrel{\alpha}{\longrightarrow} \kn_{C/\P^2_r} \longrightarrow \ko_{C'\cap
  C''} \longrightarrow 0\,,$$
$\alpha$ being induced by \mbox{$\mbox{id}_1\otimes G+F\otimes \mbox{id}_2$},
where $F$ (resp.~$G$) denotes a (local) equation of $C'$ (resp.~$C''$). Since
$C'$ and $C''$ intersect only in nodes, $\alpha$ maps precisely
\mbox{$\kn_{C'/\P^2_r}^{\;\emptyset}\oplus \kn_{C''/\P^2_r}^{\;\emptyset}$} to
\mbox{$\kn_{C/\P^2_r}^{\;\emptyset}$}, and the statement of the theorem follows
immediately (by induction) from the vanishing statement in the irreducible
case.

First, we have to make some easy considerations about exceptional curves:

\begin{lemma} \label{1.3.1}
  Let \mbox{$H\in |hE_0-\sum_{i=1}^r h_iE_i|$} be an irreducible curve, then
$$ \frac{1}{h^2} \sum_{i=1}^r h_i^2\; \leq\;1+\frac{1}{h}\;\leq 2\:.$$
\end{lemma}

\begin{proof}
The r blown--up points \mbox{$p_1,\ldots,p_r$} are chosen generically. Hence,
the existence of an
irreducible curve \mbox{$H \in |hE_0-\sum_{i=1}^r h_iE_i|$}, that is of an
irreducible curve \mbox{$\bar{H} \subset \P^2$} passing through the $p_i$ with
multiplicity $h_i$ (\mbox{$1\leq i\leq r$}), implies for an additional point
$p\:\!'\!\!_{\nu}\not\in H$, close to $p_{\nu}$ with \mbox{$h_{\nu }\geq 1$},
the existence of a
curve  \mbox{$ \bar{H'} \subset \P^2 $} passing through $p_i$ with
multiplicity $h_i$ (\mbox{$i\not= \nu$}), through $p_{\nu}$ with multiplicity
$h_{\nu}-1$ and through the additional point $p\:\!'\!\!_{\nu}$.\\
We can assume \mbox{$\,h\geq h_1\geq h_2,\ldots \geq h_r>0\,$} and obtain by
B\'ezout's theorem:
$$h_1^2+\ldots +h_{r-1}^2+h_r(h_r-1) \leq h^2\:.$$
The above statement follows immediately.
\end{proof}

\begin{remark} \label{1.3.2}
 We call an irreducible curve \mbox{$H \in |hE_0-\sum_{i=1}^r h_iE_i|$} an {\em
   exceptional curve}, if  $$\:\sum_{i=1}^r h_i^2 > h^2\:.$$ Applying
 B\'ezout's
 theorem, it is clear that for fixed data $h$, $h_i$ (\mbox{$1\leq i\leq r$})
 there is  at most one such exceptional curve $H$; hence for fixed degree $h$
 there are only finitely many exceptional curves. For example, for $h=1$ the
 exceptional curves are just the lines connecting two of the blown--up points.
\end{remark}

We divide \mbox{Sing$\,C=\Sigma_1 \cup \Sigma_2 $} where $\Sigma_2 $ denotes
the set
of all nodes lying on the ex\-cep\-tional divisors $E_i$ \mbox{$(1\leq i\leq
  r)$}. Let
\mbox{$ H\in |hE_0-\sum_{i=1}^r h_iE_i|$} be a (reduced) curve of minimal
degree passing through $\Sigma_1 $. Such a curve exists (at least) for each $h$
fulfilling
$ h(h+3)/2 \geq k $,
hence we can suppose \mbox{$h\leq [\sqrt{2k} ]\,$}. Moreover define
$$ \tilde{H}:=H \cup E_1\cup \ldots \cup E_r\:\in \:|hE_0-\sum_{i=1}^r
(h_i-1)E_i| $$
and let \mbox{$\tilde{C} \sim dE_0-\sum_{i=1}^r \tilde{d}_iE_i$} with
\mbox{$\tilde{d}_i :=
\mbox{max}\,\{d_i,\,h_i+[\frac{d_i}{2}]-1\}$}. Applying Proposition \ref{1.2.1}
we have to check three conditions:
\begin{enumerate}
\item By the Hirschowitz-Criterion (\ref{Hirschowitz})
  $\:H^1(\P_r^2,\ko_{\P^2_r}(C))$ vanishes, because (\ref{Bed2}) implies
  (\ref{1.2.3}).
\item The same criterion gives the vanishing of
  $\,H^1(\P^2_r,\ko_{\P^2_r}(\tilde{C}-\tilde{H}))$, because
  \begin{eqnarray*}
    \sqrt{\frac{(d-h+3)^2}{4}} & \geq & \frac{d+3-[\sqrt{2k}]}{2}
    \:\,\stackrel{\mbox{\footnotesize(\ref{Bed2})} }{>}
    \:\,\sqrt{1+\sum\nolimits_{i=1}^r
      \frac{(d_i+1)(d_i+2)}{2}}\\
    & \geq & \sqrt{1+\sum\nolimits_{i=1}^r
      \frac{(\tilde{d}_i-h_i+1)(\tilde{d}_i-h_i+2)}{2}}\:.\\
  \end{eqnarray*}
\item Applying Corollary \ref{1.2.6}, we have to check condition (\ref{1.2.7})
  for
  each irreducible component \mbox{$H_{\nu}\in |h^{(\nu)}-\sum_{i=1}^r h^{(\nu
    )}_iE_i|\,$}, \mbox{$1\leq \nu \leq s$}, of $H$ and each exceptional
divisor
$E_i$, \mbox{$1\leq i\leq r\:$}:
  \begin{eqnarray*}
  E_i\cdot (\tilde{C}-\tilde{H}-K_{\P_r^2} )-\#(\mbox{Sing}\,C \cap
      E_i) &\stackrel{\mbox{\footnotesize B\'ezout}}{\geq} & \tilde{d}_i -
      (h_i-1)+1-\left[\frac{d_i}{2}\right]\\
     & > & 0
  \end{eqnarray*}
  \begin{eqnarray*}
  \lefteqn{\hspace*{-0.55cm}H_{\nu}\cdot (\tilde{C}-\tilde{H} -K_{\P^2_r}) -
    \#(\mbox{Sing}\,C
   \cap H_{\nu} ) } \hspace{-0.8cm}\\
     & \stackrel{\mbox{\footnotesize B\'ezout}}{\geq} &
   \!\!\!h^{(\nu)}(d\!-\!h\!+\!3)-\! \sum_{i=1}^r
   h_i^{(\nu)}(\tilde{d}_i\!-\!(h_i\!-\!1)\!+\!1)-\!\left[\frac{h^{(\nu)}d\!-\!
   \sum_{i=1}^r h_i^{(\nu)} d_i}{2}\right]\\
       & \stackrel{\mbox{\footnotesize Cauchy}}{\geq} &
      \!\!\!h^{(\nu)}\,\left(\frac{d}{2}+3-[\sqrt{2k}]\right)-\sqrt{
        \sum\nolimits_{i=1}^r
        \big(h_i^{(\nu)}\big)^2} \:\sqrt{\sum\nolimits_{i=1}^r
        \frac{(d_i+2)^2}{4}}\\
      & \stackrel{\mbox{\footnotesize (\ref{1.3.1})}}{\geq} &
   \!\!\!h^{(\nu)}\,\left(\frac{d}{2}+3-[\sqrt{2k}]-\frac{\sqrt{2}}{2}\,\sqrt{
        \sum\nolimits_{i=1}^r (d_i+2)^2}\right)
      \;\;\stackrel{\mbox{\footnotesize (\ref{Bed1})}}{>}  \;\; 0 \;\;\qed
     \end{eqnarray*}
  \end{enumerate}

\section{Irreducibility}

\subsection{Formulation of the result}
\setcounter{equation}{0}

For $\P^2_1$, the projective plane blown up at one point $p_1$, Z.~Ran shows in
\cite{Ran} that the variety of all irreducible nodal curves
\mbox{$C\in|dE_0-d_1E_1|$}
having exactly k nodes, none of them lying on the exceptional divisor $E_1$, is
irreducible. Using the smoothness of $V_{irr}(d;d_1;k)$, one can easily deduce
its irreducibility. The aim of this section is to prove the following
irreducibility criterion for \Virr\ (\mbox{$r\geq2$}):

\begin{theorem}[Irreducibility Theorem] \label{2.1.1}
Let \mbox{$r\geq 2$}, and let the positive integers \mbox{$d;\,d_1,\ldots,d_r$}
satisfy the two (irreducibility) conditions
\renewcommand{\arraystretch}{1.5}
\begin{eqnarray}
\left[ \; \sqrt{2k}\;\, \right] & < &
\frac{d}{4}+1-\frac{1}{4}\sqrt{\sum_{i=1}^r d_i^2} \label{2.1.2} \\
\left[ \; \sqrt{2k}\;\, \right] & < &
\frac{d}{2}+1-\frac{\sqrt{2}}{2}\sqrt{\sum_{i=1}^r (d_i+2)^2}\:. \label{2.1.3}
\end{eqnarray}
Then \Virr\  is (smooth and) irreducible.

Let \mbox{$C \subset \P^2_r$} be a reduced nodal curve, \mbox{$C\in
    |D|$}. If for each irreducible component \mbox{$C_{\nu} \in
    |d^{(\nu) }E_0-\sum_{i=1}^{r} d^{(\nu)}_iE_i|\,$}, \mbox{$1\leq \nu\leq
n$},
of $C$ (having precisely $k^{(\nu)}$ nodes) the variety \mbox{$V_{irr}(d^{(\nu)
};\,d^{(\nu)}_1,\ldots,d^{(\nu)}_r;\,k^{(\nu)})$} is smooth and irreducible,
then \mbox{$V(|D|;C)$} is irreducible.
\end{theorem}

The main idea of our proof is as follows. We show that for an irreducible
curve \mbox{$C\in|dE_0-\sum_{i=1}^{r} d_iE_i|\,$} in an open dense subset
\mbox{$\tilde{V} \subset$ \Virr} the cohomology group
$H^1(C,\kn_{C/\P^2_r}^{\mbox{\scriptsize \it Sing$\,$C}})$ vanishes
(cf.~Section \ref{1.2}), especially, that the conditions imposed by fixing the
$k$ singular points are independent. It follows that the restricted morphism
$$
\begin{array}{cccl}
 \pi_{\tilde{V}}: & \tilde{V} & \longrightarrow & \mbox{Sym}^k(\P^2_r)\\
 & C & \longmapsto & \mbox{Sing}\, C \\
\end{array}
$$
is dominant, its fibres are all equidimensional and irreducible as open subsets
of the linear system
\mbox{$H^0 ( \P^2_r,\mbox{Ker}\,(\ko(dE_0\!-\!\sum
d_iE_i)\rightarrow
  \kt^1_{C,\mbox{\scriptsize \it Sing$\,$C}})\:\!)$}. Hence
$\tilde{V}$ is irreducible, which implies the irreducibility of \Virr.
The second statement is a consequence of the fact that for a fixed
reduced nodal curve $C$ as above, each generic member $\tilde{C}$ of a
component of
\mbox{$V(|D|;C)$} decomposes into components $\tilde{C}_{\nu }$, \mbox{$1\leq
\nu \leq
n$}, which are generic elements of \mbox{$V_{irr}(d^{(\nu)};\,d^{(\nu)
}_1,\ldots,d^{(\nu)}_r;\,k^{(\nu)})$}. Hence there is a well--defined dominant
morphism
$$\prod_{\nu =1}^n U_{\nu} \longrightarrow V(|D|;C)$$
where $U_{\nu}$ is open dense in \mbox{$V_{irr}(d^{(\nu)};\,d^{(\nu)
}_1,\ldots,d^{(\nu)}_r;\,k^{(\nu)})$}.

\subsection{Proof of the Irreducibility Theorem}
\setcounter{equation}{0}
 We start the proof defining the subset $\tilde{V}$ of
 \Virr\ as the set of all irreducible curves $C$ in \Virr\ having
  the subsequent properties:

\begin{enumerate}
\itemsep0.1cm
\item[(a)] \mbox{$\mbox{Sing}\,C \cap E_i =\emptyset \:$} for
  \mbox{$i=1,\ldots ,r.$}
\item[(b)] If $E$ is an exceptional curve of
  degree $e\leq 2k\,$ then \mbox{$\,\mbox{Sing}\,C \cap E =\emptyset $}.
\item[(c)] The $k$ nodes of $C$ are in general position, that is, if
  \mbox{$H\in
  |hE_0-\sum_{i=1}^r h_iE_i|\,$} is a curve containing Sing$\,C$,
\mbox{$H_{\nu}\in
  |h^{(\nu)}E_0-\sum_{i=1}^r h_i^{(\nu)}E_i|\,$} is an irreducible component of
$H$, then
  \mbox{$h^{(\nu)}(h^{(\nu)}\!+3)/2\geq k_{\nu} := \#(\mbox{Sing}\,C\cap
  H_{\nu})\:.$}
\end{enumerate}

\begin{remark} \label{2.2.1}
Condition (c) implies the existence of an irreducible curve among all
  curves \mbox{$H\in |hE_0-\sum_{i=1}^r h_iE_i|\,$} of minimal degree
  containing $\mbox{Sing}\,C$:\\
 Assume a curve \mbox{$\,H=H_1\cup H_2\cup \tilde{H}\,$} of minimal degree $h$
 decomposes \mbox{($h^{(1)}\!\leq h^{(2)}$)}, then we know that \mbox{$H_1\cup
   H_2$} contains at most
 \begin{eqnarray*}
   \lefteqn{\frac{h^{(1)}(h^{(1)}\!+3)}{2} + \frac{h^{(2)}(h^{(2)}\!+3)}{2} =
     }\hspace{0.5cm}\\
  & &\frac{(h^{(1)}+h^{(2)}\!-1)(h^{(1)}+h^{(2)}\!+2)}{2}-
\frac{h^{(1)}(h^{(2)}\!-2)+h^{(2)}(h^{(1)}\!-2)}{2}+1
 \end{eqnarray*}
nodes of $C$. The degree of $H$ being minimal, we conclude that either
\mbox{$h^{(1)}\!=1$} or  \mbox{$h^{(1)}\!=h^{(2)}\!=2$}. But in
these cases,
using the obvious constructions and Bertini's theorem, we can show the
existence of an irreducible curve of degree \mbox{$h^{(1)}\!+h^{(2)}$}, which
contains the nodes of $C$ lying on \mbox{$H_1\cup H_2$}.
\end{remark}

\begin{lemma} \label{2.2.2}
\mbox{$\,\tilde{V}\subset$ \Virr} is an open dense subset.
\end{lemma}

\begin{proof}
The openess being obvious, it is enough to show that there are no obstructions
for (locally) moving singular points of \mbox{$C\in$ \Virr} in a prescribed
position (such that the conditions (a)--(c) are satisfied). Again, we divide
\mbox{Sing$\,C=\Sigma_1\cup \Sigma_2$} where $\Sigma_2$ denotes the set of all
nodes
lying on the exceptional divisors $E_i$ \mbox{$(1\leq i\leq r$)} and start
moving nodes away from the exceptional divisors (such that
finally \mbox{$\Sigma_2=\emptyset$)}:\\
Let $z$ be a node of $C$ on $E_{i_0}$, we have to show that for
\mbox{$\Sigma_0:=\mbox{Sing}\,C-\{z\}$} the cohomology group
$H^1(C,\kn_{C/\P^2_r}^{\Sigma_0})$ vanishes: indeed, from the commutative
diagram
$$
\renewcommand{\arraystretch}{1.0}
\begin{array}{ccccccccc}
&&&& 0 \\
&&&& \downarrow \\
&&&& \ko_{\P^2_r} \\
&&&& \downarrow \\
&&&& \ko_{\P^2_r}(C) & \longrightarrow & \kt ^1_{C,\Sigma_0} & \longrightarrow
& 0  \\
&&&& \downarrow & & \| \\
0 & \longrightarrow & \kn_{C/\P^2_r}^{\Sigma_0} & \longrightarrow &
\kn_{C/\P^2_r} & \longrightarrow & \kt ^1_{C,\Sigma_0} & \longrightarrow &
0\:\, \\
&&&& \downarrow \\
&&&& 0 \\
\end{array}
$$
we can conclude the required surjectivity of
$ H^0(\P^2_r,\ko_{\P^2_r}(C))\rightarrow H^0(\P^2_r,\kt
^1_{C,\Sigma_0})\:. $

 As above, we denote by \mbox{$H\in
|hE_0-\sum_{i=1}^r h_iE_i|\,$} a curve of minimal degree \mbox{$h\leq
  [\sqrt{2k}\,]\,$}
passing through $\Sigma_1$. There are two cases to consider:

{\em Case 1} : $\,\{z\}\not\subset H$

Define $L$ to be the strict transform of a straight line
\mbox{$\overline{L}\subset \P^2$} through
$p_{i_0}$ with tangent direction corresponding to $z$ and denote
\mbox{$J:=\left\{\:j\:|\:\{p_j\}\subset \overline{L}\,\right\}\supset
\left\{i_0\right\}\:\:$} (the genericity of the blown--up points implies that
\mbox{$\#J\leq 2$}). Consider the curve \mbox{$\tilde{H}:=H\cup E_1\cup \ldots
  \cup E_r\cup
L \,\in\,|(h+1)E_0-\sum_{i=1}^r (\tilde{h}_i-1)E_i|$}, where
$$
\renewcommand{\arraystretch}{1.1}
\tilde{h}_i:=\bigg\{
 \begin{array}{c}
 h_i+1\\
 h_i
 \end{array}
 \:\mbox{for}\:
 \begin{array}{c}
 i\in J\\
 i\notin J\:.
 \end{array}
$$
Moreover, let \mbox{$\,\tilde{C}\sim dE_0-\sum_{i=1}^r \tilde{d}_iE_i\,$} with
$$\tilde{d}_i:=\mbox{max}\,\left\{d_i,\tilde{h}_i+\left[\frac{d_i}{2}\right]
- -1\right\}\:.$$
According to Proposition \ref{1.2.1} we have three conditions, sufficient
for the vanishing of $H^1(C,\kn_{C/\P^2_r}^{\Sigma_0})\,$:
as above, the Hirschowitz--Criterion (\ref{Hirschowitz}) together with
(\ref{2.1.3}) guarantees the \mbox{vanishing} property (A). (B) follows in the
same manner, because
\begin{eqnarray*}
\sqrt{\frac{(d\!-\!h\!+\!2)^2}{4} } & \geq &
\frac{d\!+\!2\!-\!\left[\sqrt{2k}\right]}{2} \:\:
 \stackrel{\mbox{\footnotesize (\ref{2.1.3})}}{>} \:\:
 \sqrt{\sum\nolimits_{i=1}^r \frac{(d_i\!+\!2)^2}{2}} \\
&\geq &\sqrt{1+\sum\nolimits_{i=1}^r
 \frac{(\tilde{d}_i\!-\!\tilde{h}_i\!+\!1)(\tilde{d}_i\!-\!\tilde{h}_i
\!+\!2)}{2} }\:.
\end{eqnarray*}
Finally, property (C) is an immediate consequence of Corollary \ref{1.2.6},
knowing
that for each exceptional divisor $E_i$ \mbox{$(1\leq i\leq r)$} we have
$$
 E_i\cdot (\tilde{C}-\tilde{H}-K_{\P_r^2} )-\#(\mbox{Sing}\,C \cap
      E_i) \stackrel{\mbox{\footnotesize B\'ezout}}{\geq}  \tilde{d}_i -
      (\tilde{h}_i-1)+1-\left[\frac{d_i}{2}\right]
      \;\:>\;\: 0
$$
and that for each irreducible component \mbox{$H_{\nu}\in
  |h^{(\nu)}-\sum_{i=1}^r h^{(\nu)}_iE_i|\,$} of
\mbox{$L\cup H$}
\begin{eqnarray*}
\lefteqn{H_{\nu}\cdot (\tilde{C}-\tilde{H} -K_{\P^2_r}) -
  \#(\mbox{Sing}\,C
  \cap H_{\nu} ) } \hspace{0.2cm}\\
 & \stackrel{\mbox{\footnotesize B\'ezout}}{\geq} &
 h^{(\nu)}\,(d\!-\!(h\!+\!1)\!+\!3)- \sum_{i=1}^r
 h_i^{(\nu)}\,(\tilde{d}_i\!-\!(\tilde{h}_i\!-\!1)\!+\!1)-
 \left[\frac{h^{(\nu)}d\!-\!\sum h_i^{(\nu)} d_i}{2}\right]\\
 & \stackrel{\mbox{\footnotesize Cauchy}}{\geq} &
 h^{(\nu)}\,\left(\frac{d}{2}+2-[\sqrt{2k}]\right)-\sqrt{ \sum\nolimits_{i=1}^r
 \big(h_i^{(\nu)}\big)^2} \:\sqrt{\sum\nolimits_{i=1}^r
 \frac{(d_i+2)^2}{4}}\\
 & \stackrel{\mbox{\footnotesize (\ref{1.3.1})}}{\geq} &
 h^{(\nu)}\,\left(\frac{d}{2}+2-[\sqrt{2k}]-\frac{\sqrt{2}}{2}\,\sqrt{
 \sum\nolimits_{i=1}^r (d_i+2)^2}\right)\:\:
 \stackrel{\mbox{\footnotesize (\ref{2.1.3})}}{>} \:\: 0\:.
\end{eqnarray*}

{\em Case 2} : $\,\{z\}\subset H$

In this case we can omit the additional component $L$ in the definition of
$\tilde{H}$ and, proceeding as in case 1, we obtain again the vanishing of
$H^1(C,\kn_{C/\P^2_r}^{\Sigma_0})\,$.

Now, we can assume \mbox{Sing$\,C=\Sigma_1$} (that is, there are no nodes of
$C$ on the exceptional divisors $E_i$) and we go on moving, subsequently, nodes
away from exceptional curves:\\
Assume \mbox{$E\in |eE_0-\sum_{i=1}^r e_iE_i|\,$} to be an irreducible curve
satisfying \mbox{$\,\sum_{i=1}^r e_i^2>e^2\,$}, let \mbox{$\,z\in \Sigma_1$} be
a node of $C$ on $E$
and denote \mbox{$\Sigma_0:=\mbox{Sing}\,C-\{z\}\,$}. As above, we construct a
curve
$$\tilde{H}:=H\cup \tilde{L}\,\in \:|(h+1)E_0-\sum_{i=1}^r h_iE_i|$$
where \mbox{$\tilde{L}\not\subset H$} is a line in $\P^2$ containing none of
the blown--up points such that \mbox{$\tilde{L}\cap
  \,\mbox{Sing}\,C=\{z\}$}. Moreover, we consider \mbox{$\,\tilde{C}\sim
  dE_0-\sum_{i=1}^r \tilde{d}_iE_i\,$} with
$$\,\tilde{d}_i:=\mbox{max}\,\left\{d_i,h_i+\left[\frac{d_i}{2}
  \right]\right\}\,,$$ and the above reasoning gives again
\mbox{$H^1(C,\kn_{C/\P^2_r}^{\Sigma_0})=0\,$}.

By Remark \ref{1.3.2}, we can end up with a curve \mbox{$C\in$
\Virr} close to the original one having properties (a) and (b). It
remains to move the nodes in general position to obtain a curve \mbox{$C\in
\tilde{V}\,$}:

Again, we choose a (not necessarily irreducible) curve \mbox{$H\in
  |hE_0\!-\!\sum_{i=1}^r h_iE_i|\,$} of minimal degree
containing \mbox{$\Sigma_1=\mbox{Sing}\,C\,$}. Assume $H$ decomposes into
irreducible components \mbox{$H_{\nu}\in |h^{(\nu)}\!-\!\sum_{i=1}^r
h^{(\nu)}_iE_i|\,$ $(1\leq \nu \leq s)$} and assume that there are more than
\mbox{$\,h^{(\nu)}(h^{(\nu)}\!+3)/2+1\,$} nodes on one component $H_{\nu}$. We
show that we
can move \mbox{$\,h^{(\nu)}(h^{(\nu)}\!+3)/2+1\,$} of them in general
position:

Define \mbox{$\tilde{H}:=H\cup G$} where \mbox{$G\in |gE_0-\sum_{i=1}^r
  g_iE_i|$}  is a curve
not containing $H_{\nu}$ through the selected
\mbox{$\,h^{(\nu)}(h^{(\nu)}\!+3)/2+1\,$} nodes. Such a curve exists for each
$g$ satisfying
$$ \frac{g(g+3)}{2} -\left(\frac{h^{(\nu)}(h^{(\nu)}\!+3)}{2}+1\right)\:>\:
\frac{(g-h^{(\nu)})(g-h^{(\nu)}+3)}{2} $$
(the right--hand side is the dimension of the linear system of curves
containing $H_{\nu}$). Hence, we can suppose \mbox{$g=h^{(\nu)}\!+2\,$}, and
\mbox{$\,h^{(\nu)}(h^{(\nu)}\!+3)/2\,+1\leq k\,$} \mbox{implies} \mbox{
$h^{(\nu)}\leq [\sqrt{2k}]-1\,$}. Proceeding as before, we have to prove the
vanishing of \mbox{$H^1(\P^2_r,\ko_{\P^2_r}(C-\tilde{H}))$}
respectively $H^1(\tilde{H},\tilde{\kn}_{\tilde{H}/\P^2_r}^{\Sigma_0})\,$. The
first
is an immediate consequence of the Hirschowitz--Criterion (\ref{Hirschowitz}),
because
\arraycolsep0.0cm
\begin{eqnarray*}
\sqrt{\frac{(d\!-\!h\!-\!g\!+\!3)^2}{4}} & \geq &\frac{d\!-\!2[\sqrt{2k}]\!+
  \!2}{2} \\
& \stackrel{\mbox{\footnotesize (\ref{2.1.3})}}{>} &
\sqrt{\sum\nolimits_{i=1}^r
  \!\frac{(d_i\!+\!2)^2}{2}}\geq \sqrt{1\!+\!\sum\nolimits_{i=1}^r
  \!\frac{(d_i\!-\!h_i\!-\!g_i)(d_i\!-\!h_i\!-\!g_i\!+\!1)}{2}}
\end{eqnarray*}
while the second results from Corollary \ref{1.2.6}, knowing that $H_{\nu}$
$(1\leq \nu \leq s)$ and the components of $G$ are no exceptional curves:
\begin{eqnarray*}
\lefteqn{H_{\nu}\cdot (C-\tilde{H} -K_{\P^2_r}) - \#(\mbox{Sing}\,C \cap
  H_{\nu} )}\hspace{0.6cm}\\
 & \stackrel{\mbox{\footnotesize B\'ezout}}{\geq} &
 h^{(\nu)}\,(d\!-\!h\!-\!g\!+\!3)- \sum_{i=1}^r
 h_i^{(\nu)}\,(d_i\!-\!h_i\!-\!g_i\!+\!1)-
 \left[\frac{h^{(\nu)}d\!-\!\sum_{i=1}^r h_i^{(\nu)} d_i}{2}\right]\\
 & \stackrel{\mbox{\footnotesize Cauchy}}{\geq} &
 h^{(\nu)}\,\left(\frac{d}{2}+2-2[\sqrt{2k}]-\frac{\sqrt{\sum\nolimits_{i=1}^r
 d_i^2}}{2}\right)\:\:
  \stackrel{\mbox{\footnotesize (\ref{2.1.2})}}{>} \:\: 0
\end{eqnarray*}
and the same holds for each irreducible component $G_{\nu}$ of G in place of
$H_{\nu}$. \qed\\
\renewcommand{\qed}{}\end{proof}

\begin{lemma} \label{2.2.3}
Let $\,C\in \tilde{V}$, then $H^1(C,\kn^{\mbox{\scriptsize \it
    Sing$\,$C}}_{C/\P^2_r})=0\,$.
\end{lemma}

\begin{proof}
By Remark \ref{2.2.1}, we can choose an irreducible curve
\mbox{$H\in |hE_0-\sum_{i=1}^r h_iE_i|\,$} of degree \mbox{$h=[\sqrt{2k}]$}
through all nodes of
$C$. Moreover, there is a curve \mbox{$G\not\supset H$}, \mbox{$G\in
  |gE_0-\sum_{i=1}^r
g_iE_i|$}, such that \mbox{Sing$\,C\subset G\cap H$}, for each $g$ satisfying
$$ \frac{g(g+3)}{2}-k\:>\:\frac{(g-h)(g-h+3)}{2}\:\:, $$
hence, especially for \mbox{$g=h+1\,$}. Let \mbox{$\tilde{H}:=H\cup
  G$}. Applying proposition \ref{1.2.1} as before, we conclude the vanishing
of $H^1(C,\kn^{\mbox{\scriptsize\it Sing$\,$C}}_{C/\P^2_r})\,$. Indeed, the
above inequalities hold again, since  neither $H$ nor the irreducible
components $G_{\nu}$ of $G$ are exceptional curves.
\end{proof}

\section{Existence}

\subsection{Formulation of the result}
\setcounter{equation}{0}
We shall treat
the problem of the existence of nodal curves in $\P^2_r$ for \mbox{$r\leq 9$}
and \mbox{$r\geq 10$}
separately. If \mbox{$r\leq 9$}, Theorems \ref{4.1A} and \ref{4.1B} will
give the complete answer,
while for \mbox{$r\geq 10$}, we obtain an asymptotically nearly optimal
sufficient criterion (Theorem \ref{3.1.1}).

\begin{theorem}[Existence Theorem A] \label{4.1A}
Let \mbox{$r=1$} then \mbox{$V_{irr}(d;d_1;k)\neq \emptyset$} if and only if
\mbox{$(d_1\leq d\!-\!1$} or \mbox{$d=d_1=1)$} and  $$0 \leq k
  \leq \frac{(d-1)(d-2)-d_1(d_1-1)}{2}.$$
Let $r=2$ and $d,d_1,d_2$ be positive integers then
  \mbox{$V_{irr}(d;d_1,d_2;k)\neq \emptyset$} if and only if $$0 \leq k
  \leq \frac{(d-1)(d-2)-d_1(d_1-1)-d_2(d_2-1)}{2}$$ and either
 \mbox{$(d_1+d_2\leq d)$} or \mbox{$(d=d_1=d_2=1)$}.
\end{theorem}

Let \mbox{$3\leq r\leq 9$}, then we define two \mbox{$(r\!+\!1)$}--tuples
\mbox{$(d;d_1,\ldots,d_r)$} and
\mbox{$(\tilde{d};\tilde{d}_1,\ldots,\tilde{d}_r)$} of non--negative integers
to be
{\em equivalent}, if there is a finite sequence of Cremona maps and
a permutation $\sigma$ transforming \mbox{$(d;d_1,\ldots,d_r)$} to
\mbox{$(\tilde{d};\tilde{d}_{\sigma(1)},\ldots,\tilde{d}_{\sigma(r)})$} . Here,
by a {\em Cremona map}, we denote a mapping
$$
\begin{array}{rccc}
\Sigma_{j,m,n} : & \Z^{r+1}&\longrightarrow & \Z^{r+1}\\
&(d;d_1,\ldots,d_r)& \mapsto & (d';d'_1,\ldots,d'_r)
\end{array}
$$
$$
\begin{array}{ll}
  \mbox{with} & d'=2d-d_j-d_m-d_n\,,\;\;\;d'_i=d_i\,\;\mbox{ for each }
  \;i\not\in\{j,m,n\},\\
& d'_j=d-d_m-d_n\,, \;\,d'_m = d-d_j-d_n\,\;\mbox{ and  }\;\,d'_n=d-d_j-d_m\,.
\end{array}
$$
Such a Cremona map corresponds to the standard Cremona transformation in
${\P}^2$ inducing the base change in Pic($\P^2_r$):
\begin{equation}
  \label{4.6}
  \left\{
\begin{array}{lcl}
E'_0 & = & 2E_0-E_j-E_m-E_n\\
E'_j & = & E_0-E_m-E_n\\
E'_m & = & E_0-E_j-E_n\\
E'_n & = & E_0-E_j-E_m\\
E'_i & = & E_i\;\;\mbox{ for each }\;i\not\in\{j,m,n\}.
\end{array}
\right.
\end{equation}
Due to the generality of the blown--up points, such a transformation maps
generic elements in \Virr\ to elements in
\mbox{$V_{irr}(d';d'_1,\ldots,d'_r;k)$} supposing
\mbox{$d,d',d_i,d'_i$} \mbox{$(1\leq i\leq r)$} to be non--negative. Since
\Virr\ is smooth, we deduce that its non--emptyness is equivalent to the
existence of a curve in \mbox{$V_{irr}(d';d'_1,\ldots,d'_r;k)$}.
An (ordered) tuple \mbox{$(d;d_1,\ldots,d_r)\in \N^{r+1}$}, \mbox{$d_1\geq
  d_2\geq \ldots \geq d_r$}, is called {\em minimal}, if it satisfies the
(minimality) condition
\begin{equation}
   \max_{\#\{j,m,n\}=3} (d_j\!+\!d_m\!+\!d_n)\: =\: d_1+d_2+d_3 \:\leq\:
     d. \label{minimal}
\end{equation}

\begin{theorem}[Existence Theorem B] \label{4.1B}
Let \mbox{$3\leq r\leq 9$} and positive integers \mbox{$d\geq d_1\geq
  \ldots\geq d_r$} satisfy the conditions
\begin{eqnarray}
\sum_{i=1}^r d_i & \leq & 3d-1 \label{4.2}\\
\sum_{i=1}^r d_i(d_i-1) & \leq & (d-1)(d-2) \label{4.3}
\end{eqnarray}
then \Virr\ $\neq \emptyset$ if and only if
\begin{eqnarray*}
0\;\:\leq\;\: k &\leq &\frac{(d\!-\!1)(d\!-\!2)}{2}-\sum_{i=1}^r
\frac{d_i(d_i-1)}{2}
\end{eqnarray*}
and \mbox{$(d;d_1,\ldots,d_r)$} is equivalent to a minimal tuple
\mbox{$(\tilde{d};\tilde{d}_1,\ldots,\tilde{d}_r)$} of non--negative integers
or to the tuple \mbox{$(1;1,1,0,\ldots,0)$}.\\
\end{theorem}

\begin{remark} \label{4.4} \noindent
\begin{enumerate}
\item Condition (\ref{4.2}) is necessary in the following
sense. By B\'ezout's Theorem and the generality of the blown--up points, the
only type of an irreducible curve not
satisfying (\ref{4.2}) is the smooth cubic through the 9 generic points.
\item For \mbox{$3\leq r\leq 8$} the tuple  \mbox{$(d;d_1,\ldots,d_r)$} is
  equivalent to a minimal one exactly if the following conditions are satisfied
  \begin{eqnarray*}
    d & \geq & d_1+d_2 \\
    2d & \geq &d_1+d_2+d_3+d_4+d_5 \\
    3d & \geq &2d_1+d_2+d_3+d_4+d_5+d_6+d_7 \\
    4d & \geq &2d_1+2d_2+2d_3+d_4+d_5+d_6+d_7+d_8 \\
    5d & \geq &2d_1+2d_2+2d_3+2d_4+2d_5+2d_6+d_7+d_8 \\
    6d & \geq &3d_1+2d_2+2d_3+2d_4+2d_5+2d_6+2d_7+2d_8 \\
  \end{eqnarray*}

\vspace{-0.7cm}
\item The exceptional case \mbox{$(d;d_1,\ldots,d_r)\sim (1;1,1,0,\ldots,0)$}
  corresponds exactly to the exceptional curves with data

  \begin{center}
  \begin{tabular}{lp{7.5cm}}
  $(2;1,1,1,1,1)$ & the conic through 5 of the generic points\\
  $(3;2,1,1,1,1,1,1)$ & the cubic through 7 of the generic points having a node
  at one of them\\
  $(4;2,2,2,1,1,1,1,1)$ & the quartic through 8 generic points having
  nodes at three of them\\
  $(5;2,2,2,2,2,2,1,1)$ & the quintic through all 8 generic points having nodes
  at 6 of them\\
  $(6;3,2,2,2,2,2,2,2)$ & the sixtic having nodes at 7 of the generic points
  and a triple point at the remaining one
  \end{tabular}
  \end{center}
\end{enumerate}

\end{remark}
For the proof, using Cremona transformations, we already saw that we can reduce
the existence problem to the case of minimal data.
Due to the independence of node smoothings in the case \mbox{$r\leq 9$},
it will then be enough to construct only rational curves, i.e.~with
$$ k=\frac{(d-1)(d-2)}{2}-\sum_{i=1}^r \frac{d_i(d_i-1)}{2} $$
nodes, in case (\ref{4.2}), (\ref{4.3}) and the minimality condition
(\ref{minimal}) are satisfied.

\begin{theorem}[Existence Theorem C] \label{3.1.1}
  Let \mbox{$r\geq 10$}, and let the positive integers
  \mbox{$d,d';\,d_1,\ldots,d_r$}
  satisfy \mbox{$d\geq d'$} and
\begin{eqnarray}
\frac{{d'}^2+6d'-1}{4} -\left[\frac{d'}{2}\right] & > & \sum_{i=1}^r
\frac{d_i(d_i+1)}{2}\:. \label{3.1.2}
  \end{eqnarray}
Then for any integer $k$ such that
\begin{eqnarray}
0\;\:\leq\;\: k &\leq &\frac{(d\!-\!1)(d\!-\!2)}{2}-\frac{d'(d'\!-\!1)}{2}\:,
 \label{3.1.3}
\end{eqnarray}
there exists a reduced irreducible curve $C$ in the linear system
\mbox{$|dE_0-\sum_{i=1}^r d_iE_i|$} on $\P^2_r$, having k nodes as its only
singularities, that is, \mbox{\Virr\ $\neq \emptyset$}.
\end{theorem}

We prove this in two steps: first, we shall prove the existence of a
non--singular curve in such linear systems on $\P^2_r$ by means of some
modification of the Hirschowitz--criterion (\ref{Hirschowitz}), afterwards, we
obtain the required nodal curves by a suitable deformation of the union of the
previous curve with generic straight lines.

\begin{corollary} \label{3.1.4}
 If \mbox{$r\geq 10$} and positive integers \mbox{$d;\,d_1,\ldots,d_r$} satisfy
$$ d\geq \sqrt{2}\sqrt{\sum_{i=1}^r d_i(d_i+1)}\:\:,$$
then, for any non--negative integer
$$ k\leq \frac{(d\!-\!1)(d\!-\!2)}{2}-\sum_{i=1}^r d_i(d_i+1)\:,$$
\Virr\ $\neq \emptyset$.
\end{corollary}

This easily follows from Theorem \ref{3.1.1}, because
\mbox{$\,d':=\sqrt{2}\sqrt{\smash[b]{\sum_{i=1}^r d_i(d_i+1)}}\;$} satisfies
$$ \frac{{d'}^2+6d'-1}{4}-\frac{d'}{2} \:>\: \sum_{i=1}^r
\frac{d_i(d_i+1)}{2}\:. $$

\subsection{Existence of nodal curves on $\P^2_r$, $r\leq9$} \label{4}
\setcounter{equation}{0}

As mentioned before, we can restrict to minimal tuples
\mbox{$(d;d_1,\ldots,d_r)$} and
construct only rational curves. In the case \mbox{$r\leq 8$} the statement is,
probably, known. Nevertheless, we provide here both, the proof for
\mbox{$r\leq8$}  and \mbox{$r=9$}.
 \subsubsection{Assume that \mbox{$r=1$}}
Clearly, B\'ezout's Theorem implies \mbox{$d_1\leq d\!-\!1$} for all
curves \mbox{$C\in V_{irr}(d;d_1;k)$}, which are not the strict transform of a
line in $\P^2$ through $p_1$.
We take the union of $d_1$ distinct straight lines in $\P^2$ through $p_1$ and
\mbox{$d-d_1$} more generic straight lines in $\P^2$, lift this curve to
$\P^2_1$ and get a reduced curve in the linear system \mbox{$|dE_0-d_1E_1|$}
with
$$\frac{(d\!-\!d_1)(d\!-\!d_1\!-\!1)}{2} +
d_1(d\!-\!d_1)=\frac{(d\!-\!1)(d\!-\!2)}{2}-\frac{d_1(d_1\!-\!1)}{2} +d\!-\!1$$
nodes. After smoothing $d-1$ intersection points, we obtain the desired
irreducible rational curve.
 \subsubsection{Assume that \mbox{$r=2$}}
 Again, by B\'ezout's Theorem, \mbox{$d_1+d_2\leq d$} with the exception of the
 strict transform of the line through $p_1$ and $p_2$. We consider in $\P^2$
 the union of $d_1$ distinct straight lines through $p_1$, $d_2$ distinct
 straight lines through $p_2$ and \mbox{$d-d_1-d_2$} additional generic
 straight lines. The strict transform in $\P^2_2$ of this curve is a reduced
 curve in the linear system \mbox{$|dE_0-d_1E_1-d_2E_2|$} with
\begin{eqnarray*}
\lefteqn{d_1d_2+(d_1\!+\!d_2)(d\!-\!d_1\!-\!d_2)
+\frac{(d\!-\!d_1\!-\!d_2)(d\!-\!d_1\!-\!d_2\!-\!1)}{2}}\hspace{1.5cm}\\
& = &\frac{(d\!-\!1)(d\!-\!2)}{2}-\frac{d_1(d_1\!-\!1)}{2}
- -\frac{d_2(d_2\!-\!1)}{2} +d\!-\!1
\end{eqnarray*}
nodes. Again one smooths \mbox{$d-1$} nodes to obtain an irreducible rational
curve.
 \subsubsection{Assume that \mbox{$r=3$}}
Due to the minimality condition (\ref{minimal}), we can proceed as follows: in
the plane we choose the union of $d_i$ distinct straight lines through $p_i$,
\mbox{$i=1,2,3$}, and \mbox{$d-d_1-d_2-d_3$} generic straight lines. After
smoothing \mbox{$d-1$} nodes of its strict transform in $\P^2_3$, as above, we
end up with the desired rational nodal curve.
 \subsubsection{Assume that \mbox{$r=4$}}
By induction on $d$, we shall show that the minimality condition
(\ref{minimal}) is sufficient for the existence of a rational irreducible
(nodal) curve
in \mbox{$|dE_0-\sum_{i=1}^4 d_iE_i|$}.
In case \mbox{$d\leq 4$} the only tuples which are not covered by the
(preceding) cases (\mbox{$r\leq 3$}) are \mbox{$(d;1,1,1,1)$}, \mbox{$d\in
  \{3,4\}$}, and \mbox{$(4;2,1,1,1)$}, hence the statement is trivial. If
\mbox{$d\geq 5$} we have
$$\max_{\#\{j,m,n\}=3} ((d_j\!-\!1)+(d_m\!-\!1)+(d_n\!-\!1))\leq d-2
$$
thus, by the induction assumption, there is an irreducible rational (nodal)
curve $C$ in the linear system
\mbox{$|(d\!-\!2)E_0-\sum_{i=1}^4 (d_i\!-\!1)E_i|$}. As is well--known
(cf.~e.g.\cite{Wae}), a
generic (smooth rational) curve $C'$ in the one--dimensional base--point--free
linear system \mbox{$|2E_0-\sum_{i=1}^4 E_i|$} intersects $C$
transversally. Finally, smoothing one intersection point in the
union of $C$ and $C'$ completes the induction step.
 \subsubsection{Assume that \mbox{$r=5$}}
We can proceed as in the case \mbox{$r=4$} with the only exception that in the
induction step $C$ will be an irreducible rational nodal curve in
\mbox{$|(d\!-\!2)E_0-d_1E_1-\sum_{i=2}^5 (d_i\!-\!1)E_i|$}, and $C'$ has to be
chosen as a generic smooth curve in \mbox{$|2E_0-\sum_{i=2}^5 E_i|$}. Thereby,
obviously, we have to treat the case
\mbox{$(d;d_1,\ldots,d_5)=(d;d\!-\!2,1,1,1,1)$} separately, because
there is no irreducible curve in the linear system
\mbox{$|(d\!-\!2)E_0-(d\!-\!2)E_1|$}. In this case, we have to choose $C$ as
the union of \mbox{$d-2$} generic lines through $p_1$ and to smooth
\mbox{$d-2$} intersection points of $C$ and $C'$.
 \subsubsection{Assume that \mbox{$r=6$}}
Again, we construct inductively irreducible rational (nodal) curves only
supposing that (\ref{minimal}) holds. First, we consider separately the case
\mbox{$d_1=\ldots=d_6=1$}, where in the induction step a generic line has to be
added (and one intersection point smoothed). Then,
supposing \mbox{$d_1\geq 2$}, in case \mbox{$d\leq 4$} the only
(additionally) possible tuple is \mbox{$(4;2,1,1,1,1,1)$}, where the statement
is trivial. For the induction step, we know, that
$$ \max \left\{(d_1\!-\!2)+(d_2\!-\!1)+(d_3\!-\!1),\:(d_2\!-\!1)+(d_3\!-\!1)+
  (d_4\!-\!1) \right\} \leq d-3\:. $$
Thereby, we can construct the desired curve by smoothing one intersection point
in the union of an irreducible nodal rational curve
$$C\in |(d\!-\!3)E_0-(d_1\!-\!2)E_1-\sum_{i=2}^6 (d_i\!-\!1)E_i|$$
and a generic curve in the one--dimensional (base--point--free) linear system
 $$ |3E_0-2E_1-E_2-E_3-E_4-E_5-E_6|\:.$$
 \subsubsection{Assume that \mbox{$r=7$}}
Changing \mbox{$(d;d_1,\ldots,d_r)$} to
\mbox{$(d\!-\!3;d_1\!-\!1,\ldots,d_7\!-\!1)$} leaves
the minimality condition intact. Hence, in the induction step (\mbox{$d\geq
  5$}), we consider the family $\kf$ of curves \mbox{$C\cup C'$}, where $C$ is
a generic rational nodal curve in
$$ |(d\!-\!3)E_0-\sum_{i=1}^7 (d_i\!-\!1)E_i| $$
and $C'$ is a generic rational nodal curve in \mbox{$|3E_0-E_1-\ldots
  -E_7|$}. First, we consider the
only case where $C$ cannot be supposed to be irreducible, namely
\mbox{$(d;d_1,\ldots,d_7)=(d;d\!-\!2,1,1,1,1,1,1)$}. In this situation, we
proceed as in case \mbox{$r=5$}, take $C$ as the union of \mbox{$d-3$} generic
lines through $p_1$ and smooth \mbox{$d-3$} intersection points of $C$ and
$C'$. Now, we assume $C$ to be irreducible. If
\mbox{$d=5$}, then $C$ is a generic curve in one of the following
base--point--free linear systems
$$|2E_0|\:,\;\;\;|2E_0-E_1|\:,\;\;\;|2E_0-E_1-E_2|\:,$$
whence a generic member of $\kf $ turns out to be a nodal curve.
It remains to prove that for \mbox{$d\geq 6$} a generic member $\tilde{C} $ of
the family $\kf$ is a nodal curve. Indeed, for the canonical divisor
$K_{\P^2_7}$ we have \mbox{$(K_{\P^2_7}\cdot C) < 0$} and
\mbox{$(K_{\P^2_7}\cdot C')<0$}, hence by (\cite{Nob}, Theorem 3.10)
\begin{eqnarray*}
\mbox{dim }\kf & = & \left(\frac{(d\!-\!3)d}{2}-\frac{(d\!-\!4)(d\!-\!5)}{2}
- -\sum_{i=1}^7 (d_i\!-\!1)\right) + \left(\frac{3\cdot6}{2}-1-7\right) + 1\\
  & = & \left(\frac{d(d\!+3)}{2}-\sum_{i=1}^7 \frac{d_i(d_i\!+\!1)}{2}\right)-
    \left(\frac{(d\!-\!1)(d\!-\!2)}{2}-\sum_{i=1}^7
      \frac{d_i(d_i\!-\!1)}{2}\right)\\
& = & \mbox{dim }|\tilde{C}| - \Delta
\end{eqnarray*}
where $\Delta$ denotes the ``virtual'' number of nodes of $\tilde{C}$. On the
other hand, the family of rational curves in \mbox{$|3E_0-E_1-\ldots -E_7|$}
has no base point. Hence, the only
possibility for a non--nodal singularity of $\tilde{C}$ is a tangency point of
$C$, $C'$ with smooth branches. But \mbox{$(K_{\P^2_7}\cdot C')=-2<-1$} and
applying (\cite{Nob}, Theorem 3.12), in this case we would have had
$$ \mbox{dim }\kf < \mbox{dim }|\tilde{C}| - \Delta \:.$$
Finally, smoothing one node in $\tilde{C}$ we get the desired curve.
 \subsubsection{Assume that \mbox{$r=8$}}
Once again, the case \mbox{$d\leq 4$} is trivial and we suppose \mbox{$d\geq
  5$}. In the induction step, we proceed as in the case \mbox{$r=7$}: We define
$\kf$ to be the family of curves \mbox{$C\cup C'$}, where C is a generic
rational nodal curve in
$$ |(d\!-\!3)E_0-\sum_{i=1}^8 (d_i\!-\!1)E_i| $$
and $C'$ is a generic rational nodal curve in \mbox{$|3E_0-E_1-\ldots
  -E_8|$}. In the case \mbox{$(d;d_1,\ldots,d_8)=(d;d\!-\!2,1,1,1,1,1,1,1)$},
we can repeat the above construction completely. Hence, by the induction
assumption,  we can assume $C$ to be irreducible. If \mbox{$d=5$}, as in the
case \mbox{$r=7$}, a generic member $\tilde{C}$ of $\kf $ is a nodal curve. If
\mbox{$d=6$} $C$, $C'$ can be chosen as the strict
transforms of two distinct (irreducible) plane rational cubics through $8$
generic points, hence \mbox{$C\cup C'$} is nodal. Let \mbox{$d\geq 7$}. Since
the minimality condition implies
$$ (K_{\P^2_8}\cdot C) = -3(d-3)+\sum_{i=1}^8 (d_i-1) \leq -3d +1
+\left[\frac{8}{3} d\right]<-1\:,$$
by (\cite{Nob}, Theorem 3.10), for dimension reasons the family of rational
nodal curves \mbox{$C\in |(d\!-\!3)E_0-\sum_{i=1}^8 (d_i\!-\!1)E_i| $} cannot
have base points. Hereby, again, the only possibility for a non--nodal
singularity of $\tilde{C}$ is a tangency point with smooth branches. Counting
dimensions as above, we see that such a singularity may not occur, and we
complete the construction as before.
\subsubsection{Assume that \mbox{$r=9$}}
By condition (\ref{4.2}) we have only to consider curves of degree \mbox{$d\geq
4$}. We split the induction step into three parts depending on the shape of
$d$:

If \mbox{$d=3m+2$}, \mbox{$m\geq 1$}, then the minimality condition
(\ref{minimal}) implies
$$ \sum_{i=1}^9 d_i\leq 3d-3 \;\;\;\mbox{ and }\;\;\;(d_1>d_3 \;\;\mbox{ or
  }\;\;d_1+d_2+d_3\leq d-1).$$
In each case, changing \mbox{$(d;d_1,\ldots,d_9)$} to
\mbox{$(d\!-\!1;d_1\!-\!1,d_2,\ldots,d_9)$} preserves both, (\ref{4.2}) and the
minimality condition. Hence, we can assume the existence of an irreducible
rational nodal curve $C$ in the linear system
 $$|(d\!-\!1)E_0-(d_1\!-\!1)E_1-\sum_{i=2}^9 d_iE_i|\:.$$
Adding a generic line through $p_1$ and smoothing one intersection point, we
obtain the desired curve.

If \mbox{$d=3m+1$}, \mbox{$m\geq 1$}, the only case where we have to use
another construction is
$$ (d;d_1,d_2,\ldots,d_9)=(3m+1;m+1,m,\ldots,m)\:. $$
We shall apply the following lemma, which will be proven at the end of this
section:
\setcounter{equation}{9}
\begin{lemma}
  \label{4.18}
Let $L_{ij}$ be the (unique) line in the linear system \mbox{$|E_0-E_i-E_j|$},
\mbox{$1\leq i < j \leq 9$}. For any \mbox{$m\geq 1$} there exists an
irreducible rational nodal
curve $$F_m\in |3mE_0-mE_1-\ldots -mE_8-(m\!-\!1)E_9|\:,$$
which meets every line $L_{ij}$, \mbox{$1\leq i<j\leq 9$}, transversally and
only at non--singular points.
\end{lemma}
Taking the curve $F_m$ from this lemma and applying the base change (\ref{4.6})
in Pic($\P^2_9$) with \mbox{$(j,m,n)=(1,2,9)$}, one easily sees that $F_m$
belongs to the linear system
$$|(3m+1)E'_0-(m+1)E'_1-(m+1)E'_2-mE_3-\ldots-mE_8-mE'_9|\:,$$
and is transversal to \mbox{$E'_2=L_{19}$}. Hence, we get the desired
irreducible rational nodal curve in
$$|(3m+1)E'_0-(m+1)E'_1-mE'_2-mE_3-\ldots-mE_8-mE'_9|\:,$$
by smoothing one intersection point of the nodal curve \mbox{$F_m\cup E'_2$}.

If \mbox{$d=3m$}, then we have to consider three possibilities
\begin{itemize}
\item \mbox{$d_1+d_2+d_3\leq d-1$}
\item \mbox{$d_1>d_3$} and \mbox{$d_1+\ldots +d_9\leq 3d-3$}
\item $(3m;m,m,m,d_4,\ldots,d_9)$, where \mbox{$m\geq d_4\geq \ldots \geq
    d_9$}, \mbox{$d_9\leq m-1$}
\end{itemize}
Clearly, to get the latter curves, it is enough to prove Lemma
\ref{4.18}. In the first two cases we can proceed as in the situation
\mbox{$d=3m+2$}.

\begin{proof}[Proof of Lemma \ref{4.18}]
  We divide our reasoning into several steps.

{\em Step 1.} First, we shall show the following. There is a hypersurface
\mbox{$S\subset (\P^2)^9$} such that on the surface $\P^2(p_1,\ldots,p_9)$
(which is the plane blown up at \mbox{$p_1,\ldots,p_9$}) the linear system
\begin{equation}
  \label{4.19}
|3mE_0-(m\!+\!1)E_1-mE_2-\ldots-mE_8-(m\!-\!1)E_9|
\end{equation}
is non--empty if and only if \mbox{$(p_1,\ldots,p_9)\in S$}.\\
Indeed, applying successively the base change (\ref{4.6}) with
\begin{equation}
  \label{4.20}
  (j,m,n)=(1,2,3),\;(4,5,6),\;(7,8,9),
\end{equation}
respectively, we transform the system (\ref{4.19}) into the system
$$|3(m\!-\!1)E'_0-mE'_1-(m\!-\!1)E'_2-\ldots-(m\!-\!1)E'_8-(m\!-\!2)E'_9|\:.$$
After $m\!-\!1$ such steps we end up with a system of type
$$|3E'_0-2E'_1-E'_2-\ldots-E'_8|$$
whose non--emptiness, evidently, imposes one condition on the nine blown--up
points. Namely, in the latter representation the condition means that the
blown--up  points \mbox{$p'_2,\ldots,p'_8$} are distinct points on a plane
cubic $C_3$
with a singularity at $p'_1$. On the other hand, choosing an irreducible cubic
$C_3$ with a singularity at a point $p'_1$, arbitrary distinct points
\mbox{$p'_2,\ldots,p'_8$} on $C_3$, and \mbox{$p'_9\not\in C_3$}, the inverse
process gives us an irreducible curve in the system (\ref{4.19}), which is
unique by B\'ezout's Theorem.

 {\em Step 2.} In the previous notation, let us specialize the points
 $p'_1,p'_4,p'_9$ on a line \mbox{$L'\subset \P^2$}. Applying the base change
 (\ref{4.20}) in the inverse order to \mbox{$\P^2(p'_1,\ldots, p'_9)$} and
 blowing down the new exceptional curves, we obtain an irreducible plane sextic
 curve with a triple point $p''_1$, double points \mbox{$p''_2,\ldots,p''_8$}
 and a non--singular point $p''_9$. Since the applied operation is a
 composition of three Cremona transformations of $\P^2$ with the fundamental
 points
 \mbox{$(p'_7,p'_8,p'_9)$},\mbox{$(p'_4,p'_5,p'_6)$},\mbox{$(p'_1,p'_2,p'_3)$},
 respectively, the
 points $p''_1,p''_4,p''_9$ lie on a straight line \mbox{$L''\subset \P^2$},
 which corresponds to the strict transform of $L'$. We continue in such a
 manner until we get an irreducible plane curve $C_{3m}(\underline{p})$ of
 degree $3m$ with an
 \mbox{$(m\!+\!1)$}--fold point $p_1$, $m$--fold points \mbox{$p_2,\ldots,p_8$}
 and an \mbox{$(m\!-\!1)$}--fold point $p_9$, such that $p_1,p_4,p_9$ lie on a
 straight line $L$. By B\'ezout's Theorem, $L$ meets $C_{3m}(\underline{p} )$
 transversally at $p_1,p_4,p_9$. Now, fixing $p'_i$, \mbox{$i\neq 4$}, we vary
 the point $p'_4$ along $C_3$. The above construction will give us a
 one--parametric family of sets \mbox{$\{p_1,\ldots, p_9\}$}, thereby a
 one--parametric (continuous) family of curves
 \mbox{$C_{3m}(p_1,\ldots,p_9)$}. Generically, the line $L$ will split into
 three lines $(p_1p_4),(p_4p_9),(p_1p_9)$, which intersect
 \mbox{$C_{3m}(\underline{p} )$} transversally. \\
 Varying the numbering of \mbox{$p'_2,\ldots, p'_8$}, respectively repeating
 the previous reasoning with the initial specialization of the points
 $p'_2,p'_5,p'_8$ on a straight line, we obtain, finally, that for a generic
 element \mbox{$\underline{p} =(p_1,\ldots,p_9)\in S$}, the strict transform of
 $C_{3m}(\underline{p} )$ intersects each $L_{ij}$, \mbox{$1\leq i<j\leq 9$}
 transversally.  Since the base change simply interchanges lines $L_{ij}$ with
 exceptional divisors $E_s$, we can claim that the intersection with each
 exceptional divisor \mbox{$E_1,\ldots,E_9$} is transversal as well.

 {\em Step 3.} The previous statement means, in particular, that for a generic
 $\underline{p} \in S$ the curve \mbox{$C_{3m}(\underline{p} )\subset \P^2$}
 has  nine ordinary multiple points and
 meets every line $(p_ip_j)$ transversally. Let us denote by $\kg$ the germ at
$C_{3m}(\underline{p} )$ of the family of plane rational curves of degree $3m$,
having ordinary singular points in a neighbourhood of \mbox{$p_2,\ldots,p_9$}
with the same multiplicities \mbox{$m,\ldots,m,m-1$}, respectively, and having
a point of multiplicity at least $m$ in a neighbourhood $U$ of $p_1$. Clearly,
$\kg$ is the intersection of a germ $\Sigma_{1} $ of the equisingular stratum
in $|\ko_{\P^2} (3m)|$, corresponding to the ordinary singular points
\mbox{$p_2,\ldots,p_9$}, and a germ $\Sigma_{2} $ of the following family of
curves of degree $3m$. A curve in $\Sigma_{2} $ has a point of multiplicity at
least $m$ in a neighbourhood $U$ of $p_1$ and the sum of $\sigma$--invariants
in $U$ is equal to $\frac{1}{2} (m\!+\!1)m$. It is not difficult to see that
$\Sigma_{2} $ is the union of $m\!+\!1$ smooth germs, such that their
intersection $\Sigma_3$ consists of curves having an ordinary
\mbox{$(m\!+\!1)$}--fold
point in $U$, and a curve in \mbox{$\Sigma_{2}\setminus \Sigma_{3}$} has in $U$
an ordinary $m$--fold point and $m$ nodes (geometrically, such a deformation
looks as if one of the local branches of $C_{3m}(\underline{p} )$ at $p_1$
moves away from the multiple point). The classical smoothness criteria say that
$\Sigma_{1} $ is smooth, and
$$ \mbox{codim}_{|\ko_{\P^2}(3m)|} \Sigma_{1} \leq 7\left(\frac{(m\!+\!1)m}{2}
  - 2 \right) + \left(\frac{m(m\!-\!1)}{2} -
  2 \right) = 4m^2+3m-16\:. $$
Since, evidently,
$$ \mbox{codim}_{|\ko_{\P^2}(3m)|} \Sigma_{2} \leq \frac{(m\!+\!1)m}{2} -2+m=
\frac{m^2}{2}+\frac{3m}{2}-2\:,$$
we obtain
 $$ \mbox{codim}_{|\ko_{\P^2}(3m)|} \kg \leq
 \frac{9m^2}{2}+\frac{9m}{2}-18\:.$$
On the other hand, we have
$$ \mbox{dim}\:(\Sigma_{1}\cap \Sigma_{3})= \mbox{dim}\:\ks = 17 < 18 \leq
\mbox{dim}\:(\Sigma_{1}\cap \Sigma_{2})\:.$$
That means, there exists a rational plane curve of degree $3m$ with $8$
ordinary $m$--fold points, one ordinary \mbox{$(m\!-\!1)$}--fold point and,
additionally, $m$ nodes. Moreover,
this curve intersects transversally with a straight line through any
two of the $9$ multiple points. Blowing up these $9$ points, we get the desired
curve \mbox{$F_m \subset \P^2_9$}.
\end{proof}

\subsection{Plane curves with generic multiple points} \label{3.2}
\setcounter{equation}{0}

It will be convenient for us to deal here with plane curves having ordinary
multiple points instead of non--singular curves on the blown--up plane. For
abuse of language, we shall use the notation: given an ordered set
\mbox{$\underline{p}=\{p_1,\ldots,p_r\}$} of distinct points in $\P^2$ and an
integral vector \mbox{$\underline{d}=(d_1,\ldots,d_r)$}, by
$S_d(\underline{p},\underline{d})$ we shall denote the set of reduced
irreducible curves of degree $d$ which have ordinary singular points at
\mbox{$p_1,\ldots,p_r$} of multiplicities
\mbox{$d_1,\ldots,d_r$}, respectively, as their only singularities.

\begin{lemma}
\label{3.2.1}
Let \mbox{$p_1,\ldots,p_r$}, \mbox{$r\geq 1$}, be distinct generic points in
$\P^2$. Then, for any positive integers \mbox{$d;d_1,\ldots,d_r$} satisfying
\begin{eqnarray}
\frac{d^2+6d-1}{4}-\left[\frac{d}{2}\right] & > & \sum_{i=1}^r
\frac{d_i(d_i+1)}{2}\:, \label{3.2.2}
\end{eqnarray}
there exists a curve \mbox{$F_d\in S_d(\underline{p},\underline{d})$}.
\end{lemma}

\begin{proof}
Following \cite{Hir}, we shall prove a more general statement. Let $I$ be a
subset in \mbox{$\{1,\ldots,r\}$}, let the points $p_i$, \mbox{$i\in I$}, lie
on a straight
line $G$, and let the points $p_i$, \mbox{$i\not\in I$}, be in general position
outside $G$. If condition (\ref{3.2.2}) and
$$ \sum_{i\in I} d_i \leq d $$
hold, then there exists a curve  \mbox{$F_d\in
  S_d(\underline{p},\underline{d})$}, which is transversal to the line $G$.

We shall use induction on $d$. If \mbox{$d=2$} then (\ref{3.2.2}) reads
$2 \geq \sum_{i=1}^r d_i(d_i+1)/2$,
so the only possibilities are: \mbox{$\,(r\!=\!1,d_1\!=\!1)$} or
\mbox{$(r\!=\!2,d_1\!=\!d_2\!=\!1)$}, when the required curves do exist. Assume
\mbox{$d\geq 3$}.

{\em Step 1}: Assume that \mbox{$I=\{1,\ldots,r\}$}, \mbox{$\,\sum_{i=1}^r
d_i\leq d\,$}.

If \mbox{$r=1$} then, by (\ref{3.2.2}), \mbox{$d_1\!<d$}, hence one obtains the
equation of the desired curve in the form
$$ F(x,y)\::=\sum_{d_1\leq i+j\leq d} A_{ij}x^iy^j\:,$$
with \mbox{$p_1=(0,0)$} and generic coefficients $A_{ij}$. If \mbox{$r>1$} then
one can obtain the desired curve as a generic member of the linear family of
all curves with equations
$$ \lambda'C'\prod_{i=1}^r \prod_{j=1}^{d_i} L'\!\!_{ij}+
\lambda''C''\prod_{i=1}^r \prod_{j=1}^{d_i}
L''\!\!\!_{ij}\:,\;\;\;(\lambda',\lambda'')\in \P^1\:,
$$
where for any \mbox{$1\leq i \leq r$}, we take
distinct generic straight lines $L'\!\!_{ij}$, $L''\!\!\!_{ij}$ through $p_i$,
\mbox{$1\leq j\leq d_i$}, and $C'$, $C''$ are distinct generic curves of degree
\mbox{$d-\sum_{i=1}^r d_i$}.

{\em Step 2}: Assume that \mbox{$\,\sum_{i=1}^r
d_i> d\:$} and \mbox{$\,\frac{d+2}{2} \leq \sum_{i\in I} d_i \leq d\:$}.

Put \mbox{$\underline{\tilde{d}}:=(\tilde{d}_1,\ldots,\tilde{d}_r)$}, where
\mbox{$\tilde{d}_i=d_i-1$}, \mbox{$i\in I$}, and \mbox{$\tilde{d}_i=d_i$},
\mbox{$i\not\in I$}.  Then
\begin{eqnarray}
\sum_{i=1}^r \frac{\tilde{d}_i(\tilde{d}_i\!+\!1)}{2} & = & \sum_{i=1}^r
\frac{d_i(d_i\!+\!1)}{2} - \sum_{i\in I} d_i \:\stackrel{\mbox{\footnotesize
    (\ref{3.2.2})}}{<}\:
\frac{d^2\!+\!6d\!-\!1}{4}-\left[\frac{d}{2}\right]-\frac{d\!+\!2}{2}
\nonumber\\
 & \leq &
 \frac{(d\!-\!1)^2+6(d\!-\!1)-1}{4}-\left[\frac{d\!-\!1}{2}\right]\:,
 \label{3.2.3}
\end{eqnarray}
hence, by the induction assumption, there exists a curve \mbox{$F_{d-1}\in
S_{d-1}(\underline{p},\underline{\tilde{d}})$}, transversal to the line
$G$. Put \mbox{$\,q:=d-\sum_{i\in I}d_i\,$} and fix \mbox{$q\!+\!1$} distinct
generic points \mbox{$z_1,\ldots,z_{q+1}$} on $G$ outside $F_{d-1}$. Since
$$ \sum_{i\in I} d_i+q+1=d+1\:,$$
we have
\begin{eqnarray}
\sum_{i=1}^r \frac{d_i(d_i+1)}{2}+q+1 & < &
\frac{d^2\!+\!6d\!-\!1}{4}-\left[\frac{d}{2}\right]+(d\!+\!1\!-\!\sum_{i\in I}
d_i ) \nonumber\\
& \leq & \frac{d^2\!+\!6d\!-\!1}{4}-\left[\frac{d}{2}\right]
+\frac{d}{2}\:\:\:\leq\:\:\:\left[\frac{(d\!+\!3)^2}{4}\right]\:,\nonumber
\end{eqnarray}
and, according to the Hirschowitz--Criterion (\ref{Hirschowitz}),
$$ h^1(\P^2,\kj(d))=0 $$
where $\kj\subset \ko_{\P^2}$ is the ideal sheaf defined by
$$\kj_{p_i}=(\frak{m}_{p_i})^{d_i}\,,\;\;i=1,\ldots,r\,,\;\;\;
\kj_{z_j}=\frak{m}_{z_j}\,,\;\;j=1,\ldots,q+1\:.
$$
That means,
\begin{eqnarray*}
h^0(\P^2,\kj(d)) & = &
h^0(\P^2,\ko_{\P^2}(d))-\sum_{i=1}^r\mbox{dim}\,\ko_{\P^2,p_i}/
(\frak{m}_{p_i})^{d_i}-\sum_{j=1}^{q+1} \mbox{dim}\,\ko_{\P^2,z_j}/
\frak{m}_{z_j} \\
& = & \frac{(d\!+\!1)(d\!+\!2)}{2}\,-\,\sum_{i=1}^r
\frac{d_i(d_i\!+\!1)}{2}-q-1\:,
\end{eqnarray*}
or, in other words, the
$$ m:=\sum_{i=1}^r \frac{d_i(d_i\!+\!1)}{2}+q+1$$
linear conditions on curves of degree $d$, imposed by the multiple points
\mbox{$p_1,\ldots,p_r,z_1,\ldots,z_{q+1}$} are independent. We write these
conditions as linear equations \mbox{$\Lambda_j(F)=0$}, \mbox{$1\leq j\leq m$},
in the coefficients of a curve $F$ of degree $d$, such that
\mbox{$\Lambda_m(F)=0$} expresses the passage of $F$ through $z_{q+1}$. Due to
the above independence, we find a curve \mbox{$F\in H^0(\ko_{\P^2}(d))$}
satisfying
$$\Lambda_1(F)=\ldots=\Lambda_{m-1}(F)=0\,,\;\;\;\Lambda_m(F)=1\:.$$

Let us consider the linear family \mbox{$\,\lambda F_{d-1}G+\mu F\,$},
\mbox{$(\lambda ,\mu)\in \P^1$}. By Bertini's \mbox{Theorem} and by the
construction of $F_{d-1}$ and $F$, the
generic member $F_{\lambda ,\mu }$ of this family is irreducible and belongs to
$S_d(\underline{p},\underline{d})$. The only thing we should show, is the
transversality of $F_{\lambda ,\mu }$ and $G$. By construction, $F_{\lambda
  ,\mu }$ has multiplicities $d_i$
at $p_i$, \mbox{$i\in I$}, respectively, and contains \mbox{$q=d-\sum_{i\in I}
  d_i$} extra points \mbox{$z_1,\ldots,z_q$} on $G$, hence, clearly,
$F_{\lambda ,\mu }$ and $G$ meet transversally.

{\em Step 3}: Assume that \mbox{$\:\sum_{i=1}^r
d_i> d\:$}, \mbox{$\;\,d_1\geq \frac{d+2}{2}\;$} and \mbox{$\;\sum_{i\in I} d_i
\leq \frac{d+1}{2}\:$}.

Define \mbox{$\underline{\tilde{d}}:=(d_1-1,d_2,\ldots,d_r)$}. As is
(\ref{3.2.3}), we obtain
$$
\sum_{i=1}^r \frac{\tilde{d}_i(\tilde{d}_i\!+\!1)}{2} <
\frac{(d\!-\!1)^2+6(d\!-\!1)-1}{4}-\left[\frac{d\!-\!1}{2}\right]\:,
$$
hence, there exists a curve \mbox{$F_{d-1} \in
  S_{d-1}(\underline{p},\underline{\tilde{d}})$}, transversal to $G$.
Note that
\begin{eqnarray*}
\mbox{dim}\,S_{d-1}(\underline{p},\underline{\tilde{d}})
& \geq & \frac{(d\!-\!1)(d\!+\!2)}{2}-\frac{(d_1\!-\!1)d_1}{2}-\sum_{i=2}^r
\frac{d_i(d_i\!+\!1)}{2} \\
& > & \frac{(d\!-\!1)(d\!+\!2)}{2}-\frac{(d\!-\!1)^2+6(d\!-\!1)-1}{4}+
\left[\frac{d\!-\!1}{2}\right] \\
& = & \frac{d^2-2d+2}{4}+\left[\frac{d\!-\!1}{2}\right]\:\:\geq \:\:2\:,
\end{eqnarray*}
as \mbox{$d\geq3$}. Therefore we find a curve \mbox{$\tilde{F}_{d-1}\in
S_{d-1}(\underline{p},\underline{\tilde{d}})$}, linearly independent of
$F_{d-1}$. Consider the linear family
$$\lambda F_{d-1}L+\mu \tilde{F}_{d-1} \tilde{L}\,,\;\;\;(\lambda, \mu)\in
\P^1\,,$$
where $L$, $\tilde{L}$ are distinct generic straight lines through $p_1$. Then,
by Bertini's Theorem a generic member of this family is irreducible, belongs to
$S_d(\underline{p},\underline{d})$ and is transversal to $G$.

{\em Step 4}: Assume that $d$ is odd and \mbox{$\;\sum_{i\in I} d_i
=\frac{d+1}{2}=d_m:=\mbox{max}\,\{d_i\,|\:i\not\in I\,\}\,$}.

First, we show that $\#(I)>1$. Indeed, otherwise, we had
$$ 0\,\stackrel{\mbox{\footnotesize (\ref{3.2.2})}}{<}\,
\frac{d^2\!+\!6d\!-\!1}{4}-\frac{d\!-\!1}{2}
- -\frac{d\!+\!1}{2}
  \left(\frac{d\!+\!1}{2}+1\right)\: =\: -\frac{1}{2}\:.
$$
Choose \mbox{$j\neq k \in I$}. Since
\mbox{$d_m\!+\min\,\{d_j,\,d_k\}\geq\frac{d+2}{2}\,$} and
\mbox{$d_m\!+\max\,\{d_j,\,d_k\}\leq d\,$}, once again as in
(\ref{3.2.3}), we obtain for
$$
\begin{array}{lll}
\underline{\tilde{d}}=(\tilde{d}_1,\ldots,\tilde{d}_r)\,, &
\tilde{d}_i:=d_i\,\mbox{ for } i\not\in
\{m,j\}\,, & \tilde{d}_i:=d_i\!-\!1\,
\mbox{ for } i\in \{m,j\} \\
\underline{\tilde{d'}\!}=(\tilde{d}'\;\!\!\!_1,\ldots,
\tilde{d}'\;\!\!\!_r)\,,&
\tilde{d}'\;\!\!\!_i:=d_i\,\mbox{ for } i\not\in
\{m,k\}\,, & \tilde{d}'\;\!\!\!_i:=d_i\!-\!1\,
\mbox{ for } i\in \{m,k\}
\end{array}
$$
the existence of curves  \mbox{$\tilde{F}_{d-1}\in
S_{d-1}(\underline{p},\underline{\tilde{d}}\,)$},
\mbox{$\,\tilde{F}'\!\!_{d-1}\in
  S_{d-1}(\underline{p},\underline{\tilde{d'}\!}\:)$}, transversal to $G$. Let
$\tilde{G}$ be the straight line through $p_m$, $p_j$, and $\tilde{G'}$ be the
straight line through $p_m$, $p_k$. Consider the linear family
$$ \tilde{\lambda}\tilde{F}_{d-1}\tilde{G} + \tilde{\lambda
  '}\tilde{F}'\!\!_{d-1}\tilde{G'},\;\;\;(\tilde{\lambda},\tilde{\lambda
  '})\in \P^1\:.$$
By Bertini's Theorem, a generic member $F$ of this family has the only singular
points \mbox{$p_1,\ldots,p_r$} of multiplicities \mbox{$d_1,\ldots,d_r$},
respectively, and is transversal to $G$. Since
\mbox{$\tilde{F}_{d-1}\tilde{G}$}  has ordinary
singularities at $p_i$, \mbox{$i\not\in \{j,m\}$}, and
\mbox{$\tilde{F}'\!\!_{d-1}\tilde{G'}$}  has ordinary
singularities at $p_i$, \mbox{$i\not\in \{k,m\}$}, $F$ has ordinary
singularities at $p_i$, \mbox{$i\neq m$}. At the point $p_m$, the curve
\mbox{$\tilde{F}_{d-1}\tilde{G}$} has at most one multiple tangent, which
should coincide with $\tilde{G}$, on the other hand
\mbox{$\tilde{F}'\!\!_{d-1}\tilde{G'}$} has at most one multiple tangent, which
should coincide with \mbox{$\tilde{G'}\neq \tilde{G}$}. Therefore, $F$ has no
multiple tangent at $p_m$. Finally we have to show that $F$ is irreducible,
but this follows immediately from Bertini's Theorem: indeed, the only
possibility for $F$ to be reducible is \mbox{$F=F_{d-1}L$},
where $F_{d-1}$ is an irreducible curve of degree \mbox{$d\!-\!1$} and $L$ is a
straight line, which must vary from $\tilde{G}$ to $\tilde{G'}$ as
$(\tilde{\lambda},\tilde{\lambda'})$ runs through $\P^1$, which is impossible,
because $F_{d-1}$ must have multiplicity $d_j$ at $p_j$ and multiplicity $d_k$
at $p_k$, while $\tilde{F}_{d-1}$ has multiplicity \mbox{$d_j-1$} at $p_j$.

{\em Step 5}: Assume that \mbox{$\:\sum_{i=1}^r
d_i> d\:$}, \mbox{$\;\sum_{i\in I} d_i\leq
\min\,\{\frac{d+1}{2},\,d\!-\!d_m$\}} and \mbox{$\;d_m\leq \frac{d+1}{2}\,$}.

(Again, $d_m$ denotes \mbox{$\max\,\{d_i\,|\:i\not\in I\}$}.) In this case, we
specialize the point $p_m$ to a generic point in
\mbox{$G\setminus\{p_i\,|\:i\in I\}$}, and we end up with one of the cases
2--5.

Finally, consider for each $r$--tuple \mbox{($p_1,\ldots,p_r$)} of points in
$\P^2$ the non--empty linear system of curves with multiplicity at least $d_i$
at $p_i$ (\mbox{$1\leq i\leq r$}). Then in each of the occurring cases the
Hirschowitz criterion (\ref{Hirschowitz}) implies the non--speciality. Hence
such linear systems are equidimensional and give an irreducible variety. But
the condition to be irreducible with only given non--degenerate multiple points
is open (it is described by inequalities); hence the existence of such a curve
in a more special situation (as considered above) implies the existence in the
original situation.
\end{proof}

\subsection{Construction of nodal curves} \label{3.3}
\setcounter{equation}{0}
Given \mbox{$k,d,d',d_1,\ldots,d_r$} satisfying conditions (\ref{3.1.2}) and
(\ref{3.1.3}), we shall construct a reduced irreducible plane curve of degree
$d$ with ordinary singular points \mbox{$p_1,\ldots,p_r$} of multiplicities
\mbox{$d_1,\ldots,d_r$}, respectively, and with $k$ additional nodes as its
only singularities.

Let us fix distinct generic points \mbox{$p_1,\ldots,p_r\in \P^2$}, take an
irreducible curve
\mbox{$\Phi\in S_{d'}(\underline{p},\underline{d})$}, and put
$$F:=\Phi \prod_{i=1}^{d-d'} L_i\:,$$
where $L_1,\ldots,L_{d-d'}$ are distinct generic straight lines. This curve $F$
has ordinary singular points \mbox{$p_1,\ldots,p_r$} of multiplicities
\mbox{$d_1,\ldots,d_r$}, respectively, and \mbox{$m:=d(d-1)/2-d'(d'-1)/2$}
nodes \mbox{$z_1,\ldots,z_m$} as its only singularities. We shall show that it
is possible to smooth prescribed nodes keeping the given ordinary singularities
and the rest of nodes, and thus prove Theorem \ref{3.1.1}.

 As we have seen
above, we can deduce the required independence of the deformations of the nodes
from
\begin{eqnarray}
H^1(F,\kn^{\,'}\!\!\!_{F/\P^2})=0\:,\label{3.3.1}
\end{eqnarray}
\mbox{$ \kn^{\,'}\!\!\!_{F/\P^2}:=\mbox{Ker}\,(\kn_{F/\P^2}\rightarrow
  \kt_F)\:,$} where $\kt_F$ denotes the skyscraper sheaf on $\P^2$ concentrated
at \mbox{$p_1,\ldots,p_r$}, \mbox{$z_1,\ldots,z_m$} defined by
$$
\kt_{F,p_i}:=\ko_{\P^2,p_i}/(\frak{m}_{p_i})^{d_i}\,,
\;\:1\leq i\leq r\,,\:\:\;\;
\kt_{F,z_j}:=\ko_{\P^2,z_j}/\frak{m}_{z_j}\,,\;\:1\leq j\leq m\,.$$
We prove (\ref{3.3.1}) by induction on $d$. If \mbox{$d=d'$}, then the
vanishing of
$ H^1(F,\kn^{\,'}\!\!\!_{F/\P^2}) $ is provided by the Hirschowitz--Criterion
(\ref{Hirschowitz}), because (\ref{3.1.2}) implies
$$\left[\frac{(d'\!+3)^2}{4}\right] > \sum_{i=1}^r\frac{d_i(d_i+1)}{2}\:.$$
Assume that \mbox{$d>d'$}. Then denote
$$ \tilde{F}:=\Phi\prod_{i=1}^{d-d'-1} L_i \,,$$
that is,
  \mbox{$F=\tilde{F}L_{d-d'}$}.
 Let \mbox{$\tilde{F}\cap
L_{d-d'}:=\{z_1,\ldots z_l\}$}, then $\kt_{\tilde{F}}$ denotes the restriction
of $\kt_F$ on \mbox{$\{p_1,\ldots,p_r,z_{l+1},\ldots,z_m\}$}. Consider the
exact sequence
\renewcommand{\arraystretch}{1.0}
$$
\begin{array}{ccccccccc}
0&\longrightarrow &\ko_{\tilde{F}}(d-1) \oplus \ko_{L_{d-d'}}(1)
&\stackrel{\alpha}{\longrightarrow} & \ko_F(d)
&\longrightarrow &\ko_{\tilde{F}\cap L_{d-d'}} & \longrightarrow  & 0\\
& & \:\| & & \| \\
& & \:\:\kn_{\tilde{F}/\P^2}\oplus \ko_{\P^1}(1) & & \kn_{F/\P^2} \\
\end{array}
$$
where $\alpha$ is given by \mbox{$\,\mbox{id}_1\otimes
L_{d-d'}+\tilde{F}\otimes \mbox{id}_2\,$}. Clearly, $\alpha $ maps
\mbox{$\kn^{\,'}\!\!\!_{\tilde{F}/\P^2}\oplus \ko_{\P^1}(1)$} to
$\kn^{\,'}\!\!\!_{F/\P^2}$. Therefore, (\ref{3.3.1})
follows from the induction assumption
$$ H^1(\tilde{F},\kn^{\,'}\!\!\!_{\tilde{F}/\P^2})=0\:, $$
and we are finished.


\end{document}